\documentclass[twoside,a4paper,12pt]{article}

\usepackage{amsmath,amsthm,amssymb}
\usepackage[v2,tips]{xy}

\newcommand{\mc}[1]{\mathcal{#1}}
\newcommand{\ip}[2]{{\langle {#1} , {#2} \rangle}}
\newcommand{\aone}{\Box}
\newcommand{\atwo}{\Diamond}
\newcommand{\proten}{{\widehat{\otimes}}}
\newcommand{\noproof}{\hfill$\square$}
\newcommand{\wap}{\operatorname{WAP}}

\newcommand{\tobidual}[1]{\kappa_{#1}(#1)}

\newcommand{\topcen}[2]{{\mf Z_t^{(#1)}(#2)}}
\newcommand{\mf}[1]{\mathfrak{#1}}
\newcommand{\conv}{\operatorname{convex}}
\newcommand{\inten}{{\check{\otimes}}}
\newcommand{\bimodhomo}[2]{{_{#1}}\mc B_{#1}(#2)}

\theoremstyle{plain}
\newtheorem{proposition}{Proposition}[section]
\newtheorem{theorem}[proposition]{Theorem}
\newtheorem{corollary}[proposition]{Corollary}
\newtheorem{lemma}[proposition]{Lemma}
\theoremstyle{definition}
\newtheorem{definition}[proposition]{Definition}

\theoremstyle{remark}
\newtheorem{example}[proposition]{Example}
\newtheorem{remark}[proposition]{Remark}

\begin{document}

\title{Dual Banach algebras: representations and injectivity}
\author{Matthew Daws}
\maketitle

\begin{abstract}
We study representations of Banach algebras on reflexive Banach spaces.
Algebras which admit such representations which are bounded below seem
to be a good generalisation of Arens regular Banach algebras; this class
includes dual Banach algebras as defined by Runde, but also all group
algebras, and all discrete (weakly cancellative) semigroup algebras.
Such algebras also behave in a similar way to C$^*$- and W$^*$-algebras;
we show that interpolation space techniques can be used in the place of GNS
type arguments.  We define a notion of injectivity for dual Banach algebras,
and show that this is equivalent to Connes-amenability.  We conclude by
looking at the problem of defining a well-behaved tensor product
for dual Banach algebras.
\end{abstract}

{\small\noindent
\emph{Keywords:} dual Banach algebra, von Neumann algebra, Connes-amenability,
group algebra, unique predual}

{\small\noindent
2000 \emph{Mathematical Subject Classification:}
46B70, 46H05, 46H15, 46H99, 46M05, 47L10 (primary),
43A10, 43A20, 46A25, 46A32, 46A35, 46L10, 46L06, 46M20.}


\section{Introduction}

It has been known for some time (see \cite{Kai} and \cite{Young})
that a Banach algebra $\mc A$ which admits a faithful representation on
a reflexive Banach space has an intrinsic characterisation, namely
that the \emph{weakly almost periodic functionals}, written $\wap(\mc A')$
(see below for the definition), separate the points of $\mc A$.  Similarly,
if we wish to find an isometric representation of this kind, we need only
ask that $\wap(\mc A')$ form a norming set for $\mc A$.  We shall call
$\mc A$ a \emph{WAP-algebra} when $\mc A$ admits an isomorphic representation
on a reflexive Banach space.  Such algebras seem not
to have been studied abstractly before, but they seem to be a good
generalisation of Arens regular Banach algebras, and to form a good
framework for studying dual Banach algebras.

We follow the notation of \cite{Dales}, writing $\ip{\cdot}{\cdot}$
for the dual pairing between a Banach space $E$ and its dual, $E'$.
We write $\kappa_E:E\rightarrow E''$ for the canonical map given by
$\ip{\kappa_E(x)}{\mu} = \ip{\mu}{x}$ for $x\in E$ and $\mu\in E'$.
When $E$ is reflexive, we tend to identify $E$ with $E''$.  We write
$\mc B(E,F)$ for the space of all bounded linear operators between Banach
spaces $E$ and $F$, and we denote by $\mc F(E,F), \mc A(E,F), \mc K(E,F)$
and $\mc W(E,F)$ the subspaces of, respectively, finite-rank, approximable,
compact and weakly-compact operators (so $\mc A(E,F)$ is the operator-norm
closure of $\mc F(E,F)$ in $\mc B(E,F)$).  We write $\mc B(E)$ for
$\mc B(E,E)$, and so forth.

A \emph{dual Banach algebra} is a Banach algebra $\mc A$ such that
$\mc A = E'$, as a Banach space, for some Banach space $E$, and such
that the multiplication on $\mc A$ is separately weak$^*$-continuous.
Recall that a W$^*$-algebra is a C$^*$-algebra which is a dual
Banach algebra.  However, it is known that the multiplication (and the
involution) are automatically weak$^*$-continuous in this case.  We
use \cite{tak} as general references for C$^*$-
and W$^*$-algebras.  Dual Banach algebras were introduced in
\cite{Runde1}, but had been studied previously under different names.

For a Banach algebra $\mc A$, we turn $\mc A'$ into an $\mc A$-bimodule
in the obvious way, by setting
\[ \ip{a\cdot\mu}{b} = \ip{\mu}{ba},\quad
\ip{\mu\cdot a}{b} = \ip{\mu}{ab} \qquad (a,b\in\mc A,\mu\in\mc A'). \]
We may then check that for a dual Banach algebra $\mc A=E'$, we have that
$\kappa_E(E)$ is a submodule of $\mc A'=E''$.  We call $E$ the \emph{predual}
of $\mc A$, and write $(\mc A,E)$ if we wish to stress which predual
we are using, and often write $\mc A_*$ for $E$.  In this special case,
we shall often suppress the map $\kappa_{\mc A_*}$, and speak of $\mc A_*$
as being a subspace of $\mc A'$ (see Definition~\ref{other_defn} for
a justification of this).
We later study when such preduals are unique, both in
the isometric sense (as is well-known for W$^*$-algebras) and the
isomorphic sense, which seems more natural for Banach algebras.

When $E$ is a reflexive Banach space, the \emph{projective tensor
product} of $E$ with its dual $E'$, denoted by $E'\proten E$, is the
canonical predual for $\mc B(E)$ (see \cite{Runde1}).  This induces
a weak$^*$-topology on $\mc B(E)$.
Recall that the norm on $E'\proten E$ is $\pi(\cdot)$, defined by
\[ \pi(\tau) = \inf\Big\{ \sum_{k=1}^n \|x_k\| \|y_k\|
: \tau = \sum_{k=1}^n x_k \otimes y_k \Big\}
\qquad (\tau\in E'\otimes E). \]
We write $E'\proten E$ and not $E\proten E'$ as the former makes more
sense when $E$ is not necessarily reflexive; the two spaces are isometrically
isomorphic.  Here and elsewhere, we refer the reader to \cite{Ryan},
\cite{DF} or \cite[Chapter~VIII]{DU} for further details on tensor products
of Banach spaces.

Let $\mc A$ be a Banach algebra, and for $\mu\in\mc A'$, define $L_\mu,
R_\mu\in\mc B(\mc A,\mc A')$ by
\[ L_\mu(a) = \mu\cdot a, \quad R_\mu(a) = a\cdot\mu
\qquad (a\in\mc A). \]
Then $\mu\in\wap(\mc A')$ if and only if $L_\mu\in\mc W(\mc A,\mc A')$
(which is equivalent to $R_\mu\in\mc W(\mc A,\mc A')$).  This notation
differs from that sometimes used, but follows \cite[Section~4]{Runde2},
for example.  It may be easily checked (as we do below) that for a
dual Banach algebra $(\mc A,\mc A_*)$, we have that $\mc A_*\subseteq
\wap(\mc A')$.  It hence immediately follows from \cite{Young} that
there exists a reflexive Banach space $E$ and an isometric representation
$\pi:\mc A\rightarrow\mc B(E)$ (here \emph{representation} simply means
a homomorphism to a Banach algebra of the form $\mc B(E)$).  However, it
is not immediately apparent if such a representation need be
weak$^*$-continuous.  As usual, we may regard representations and
left-modules as interchangeable, so this question is equivalent to $E$
being a \emph{normal} in the sense of \cite{Runde2}.  We show below
that we can indeed choose $E$ to be normal (actually, our argument is very
similar to that used by Kaiser and Young, but the required machinery,
interpolation spaces, shall be needed later anyway).  This shows that
dual Banach algebras can be thought of as ``abstract'' weak$^*$-closed
subalgebras of $\mc B(E)$ for reflexive Banach spaces $E$.  This exactly
mirrors the fact that W$^*$-algebras are abstract von Neumann algebras,
that is, weak$^*$-closed (which in this context agrees with weak operator
topology closed) self-adjoint subalgebras of $\mc B(H)$ for a Hilbert space $H$.

A \emph{derivation} from a Banach algebra $\mc A$ to an $\mc A$-bimodule
$E$ is a bounded linear map $d$ such that $d(ab)=a\cdot d(b)+d(a)\cdot b$.
Fix $x\in E$, and defined $d$ by $d(a)=a\cdot x-x\cdot a$.  Then $d$ is
a derivation, called an \emph{inner derivation}.  We say that $\mc A$ is
\emph{amenable} if every derivation from $\mc A$ to a dual
$\mc A$-bimodule $E'$ is inner.  We refer the reader to \cite{RundeBook}
for details on amenability.  Similarly, Runde defines a dual Banach
algebra $(\mc A,\mc A_*)$ to be \emph{Connes-amenable} if every
weak$^*$-continuous derivation to a normal, dual $\mc A$-bimodule is inner.

For a W$^*$-algebra $\mc A$, it is this notion of amenability which
seems most natural.  One of the major achievements of C$^*$-algebra theory
has been to give equivalent natural conditions for a W$^*$-algebra to
be Connes-amenable (see \cite[Chapter~6]{RundeBook}).  One of these is
the notation of \emph{injectivity}.  We define a similar (though weaker)
notion for dual Banach algebras, and show that it is equivalent to
Connes-amenability.

We finish the paper with a study of tensor products of dual Banach
algebras.  This last section is slightly more speculative, but it
is the author's opinion that fully understanding tensor products seems
central to understanding notions of amenability: certainly the rather
well-behaved tensor products of C$^*$-algebras play a central role
in the theory of amenability for such algebras (for example, the fact
that amenability is equivalent to nuclearity).

\section{Basic properties of WAP and dual Banach algebras}

In this section, we shall study the basic properties of dual Banach
algebras, and define WAP-algebras.

Following, for example, \cite{Dales}, in dealing with Banach algebras,
we assume, by means of standard renormings, that the product is
contractive (and not merely bounded) and that a unit always has norm one.
This philosophy is compatible with dual Banach algebras:

\begin{proposition}
Let $E$ be a Banach space such that $\mc A = E'$ admits a bounded
algebra product.  Then there is an equivalent norm on $E$ such that
$\mc A$ becomes a Banach algebra.  If $\mc A$ has a unit $e_{\mc A}$,
we may choose this norm such that $\|e_{\mc A}\|=1$.
\end{proposition}
\begin{proof}
Suppose that $1 < M = \sup\{ \|ab\| : a,b\in\mc A, \|a\|=\|b\|=1\}$,
for if $M\leq 1$, we have nothing to do.  Then define
$\|\mu\|_0 = M^{-1} \|\mu\|$ for $\mu\in E$.
For $a\in\mc A$, we then have that $\|a\|_0 =
\sup\{ |\ip{a}{\mu}| : \|\mu\|\leq M \} = M\|a\|$, so
that $\|ab\|_0 = M\|ab\| \leq M^2 \|a\| \|b\| = \|a\|_0\|b\|_0$.

Now suppose that $\mc A$ has a unit $e_{\mc A}$.
Let
\[ X = \conv \{ a\cdot\mu : \|a\|_0 = \|\mu\|_0 \leq 1 \}
\subseteq E, \]
so that for $\mu\in E$, we have that $\mu = e_{\mc A}\cdot\mu
\in \|e_{\mc A}\|_0 \|\mu\|_0 X$.  Thus we can define
$\|\cdot\|_1$ on $E$ by
\[ \|\mu\|_1 = \inf\{ t>0 : \mu\in tX \} \qquad (\mu\in E), \]
and have that $\|\mu\|_1 \leq \|e_{\mc A}\|_0 \|\mu\|_0$.
Conversely, we have that $X \subseteq \{ \mu\in E :
\|\mu\|_0\leq 1\}$, so if $\|\mu\|_1=1$, then for each
$\epsilon>0$, we have that $\mu \in (1+\epsilon)X \subseteq
\{ (1+\epsilon)\lambda : \lambda\in E, \|\lambda\|_0\leq1\}$,
and so $\|\mu\|_0 \leq 1$.  Thus $\|\cdot\|_1$ is equivalent to
$\|\cdot\|_0$ and hence also equivalent to $\|\cdot\|$.

Then, for $a\in\mc A$, we have that
\begin{align*}
\|a\|_1 &= \sup\{ |\ip{a}{\mu}| : \mu\in X \}
= \sup\Big\{ \Big|\sum_{j=1}^n \ip{a}{b_j\cdot\mu_j}\Big| :
   \sum_{j=1}^n \|b_j\|_0\|\mu_j\|_0 \leq 1 \Big\} \\
&= \sup\{ |\ip{a}{b\cdot\mu}| : \|b\|_0 = \|\mu\|_0 \leq 1 \}
= \sup\{ \|ab\|_0 : \|b\|_0=1 \}.
\end{align*}
Note that $\|ab\|_0 \leq \|a\|_1 \|b\|_0$ for $a,b\in\mc A$.
We hence have that
\[ \|ab\|_1 = \sup\{ \|abc\|_0 : \|c\|_0=1 \}
\leq \|a\|_1 \sup\{ \|bc\|_0 : \|c\|_0=1 \}
= \|a\|_1 \|b\|_1, \]
and clearly $\|e_{\mc A}\|_1=1$, as required.
\end{proof}

Let $(\mc A,\mc A_*)$ be a dual Banach algebra, and let
$\mc A_*^\flat$ be the Banach space $\mc A_* \oplus \mathbb C$ with
norm
\[ \| (\mu,\alpha) \| = \max( \|\mu\|, |\alpha| )
\qquad (\mu\in\mc A_*, \alpha\in\mathbb C). \]
Then $(\mc A_*^\flat)' = \mc A\oplus\mathbb C = \mc A^\flat$ with norm
\[ \| (a,\beta)\| = \|a\|+|\beta|
\qquad (a\in\mc A, \beta\in\mathbb C). \]
We turn $\mc A^\flat$ into a Banach algebra by setting
$(a,\alpha)(b,\beta) = (ab+\beta a + \alpha b,\alpha\beta)$.
It is a simple verification that then $(\mc A^\flat,\mc A_*^\flat)$
is a dual Banach algebra.

We set $(\mc A,\mc A_*)^\sharp$ to be $(\mc A,\mc A_*)$ when
$\mc A$ is unital, and to be $(\mc A^\flat,\mc A_*^\flat)$ otherwise.
This gives us a (rather crude, it turns out) way to unitise a dual
Banach algebra.

\begin{lemma}
Let $E$ be a Banach space such that $\mc A = E'$ is a Banach algebra.
Then $(\mc A,E)$ is a dual Banach algebra if and only if $\kappa_E(E)
\subseteq \mc A'$ is a sub-$\mc A$-bimodule.
\end{lemma}
\begin{proof}
This is a routine calculation showing that the product is separately
weak$^*$-continuous if and only if $E$ is an $\mc A$-bimodule.
\end{proof}

We shall now recall the \emph{Arens products} which shall allow us
to prove some simple facts about dual Banach algebra (much in the
spirit of \cite{Palmer}).  Most of the following results
are folklore (compare, for example, with \cite[Section~1]{LL2})
but do not appear to have formally been collected together before.

For a Banach algebra $\mc A$, we turn
$\mc A'$ into a Banach $\mc A$-bimodule in the standard way (and
hence $\mc A''$ as well).  We then define bilinear maps $\mc A''\times
\mc A' \rightarrow \mc A'$ and $\mc A'\times\mc A''\rightarrow\mc A'$ by
\[ \ip{\Phi\cdot\mu}{a} = \ip{\Phi}{\mu\cdot a}, \quad
\ip{\mu\cdot\Phi}{a} = \ip{\Phi}{a\cdot\mu} \qquad
(a\in\mc A, \mu\in\mc A', \Phi\in\mc A''). \]
We then define two bilinear maps $\aone,\atwo:\mc A''\times\mc A''
\rightarrow\mc A''$ by
\[ \ip{\Phi\aone\Psi}{\mu} = \ip{\Phi}{\Psi\cdot\mu}, \quad
\ip{\Phi\atwo\Psi}{\mu} = \ip{\Psi}{\mu\cdot\Phi}
\qquad (\mu\in\mc A', \Phi,\Psi\in\mc A''). \]
We can then calculate that
\[ (\mu\cdot\Phi)\cdot\Psi = \mu\cdot(\Phi\atwo\Psi), \quad
\Phi\cdot(\Psi\cdot\mu) = (\Phi\aone\Psi)\cdot\mu
\qquad (\mu\in\mc A', \Phi,\Psi\in\mc A''), \]
from which it follows that $\aone$ and $\atwo$ are Banach-algebra
products, called the \emph{first} and \emph{second Arens products}
respectively (see \cite[Section~1.4]{PalBook} or
\cite[Theorem~2.6.15]{Dales} for further
details).  Further, we have that $\kappa_{\mc A}(a) \aone \Phi =
a\cdot\Phi = \kappa_{\mc A}(a) \atwo \Phi$ for $a\in\mc A,\Phi\in\mc A''$,
and similarly $\Phi\aone \kappa_{\mc A}(a) = \Phi\cdot a =
\Phi\atwo\kappa_{\mc A}(a)$.  When $\aone = \atwo$, we say that
$\mc A$ is \emph{Arens regular}, and in this case, we may check that
$(\mc A'',\mc A')$ becomes a dual Banach algebra.

\begin{lemma}\label{wap_char}
Let $\mc A$ be a Banach algebra, and let $\mu\in\mc A'$.  Then
$\mu\in\wap(\mc A')$ if and only if $R_\mu\in\mc W(\mc A,\mc A')$,
which is if and only if $\ip{\Phi\aone\Psi}{\mu} = \ip{\Phi\atwo\Psi}{\mu}$
for $\Phi, \Psi\in\mc A''$.
\end{lemma}
\begin{proof}
This is a simple calculation: see \cite[Proposition~3.11]{DL}.
A key tool is Gantmacher's theorem, which states that $T\in\mc W(E,F)$ if
and only if $T''(E'') \subseteq \kappa_F(F)$.
\end{proof}

We also define the two \emph{topological centres} (see \cite{LU}
or \cite{DL}) by
\begin{align*}
\topcen{1}{\mc A''} &= \{ \Phi\in\mc A'' : \Phi\aone\Psi =
   \Phi\atwo\Psi \ (\Psi\in\mc A'') \}, \\
\topcen{2}{\mc A''} &= \{ \Phi\in\mc A'' : \Psi\aone\Phi =
   \Psi\atwo\Phi \ (\Psi\in\mc A'') \}.
\end{align*}
Then the Arens products agree on either of the topological centres,
each topological centre is an algebra, and 
$\topcen{1}{\mc A''} \cap \topcen{2}{\mc A''} \supseteq
\kappa_{\mc A}(\mc A)$ is an ideal in $\mc A''$ with respect
to either Arens product.

For a Banach space $E$, a subspace $F$ of $E$, and a subspace
$G$ of $E'$, we define
\[ F^\perp = \{ \mu\in E' : \ip{\mu}{x}=0 \ (x\in F) \}, \quad
{^\perp} G = \{ x\in E : \ip{\mu}{x}=0 \ (\mu\in G) \}. \]
It is then standard that, when $F$ is closed, $F'$ is
isometrically isomorphic to $E' / F^\perp$, while $(E/F)'$
is isometrically isomorphic to $F^\perp$.  The weak$^*$-closure
of $G$ in $E'$ is $({^\perp}G)^\perp$.

\begin{proposition}\label{wap_dba}
Let $\mc A$ be a Banach algebra, and let $X \subseteq \mc A'$
be a closed submodule.  Then the following are equivalent:
\begin{enumerate}
\item the first Arens product drops to a well-defined product
   on $X' = \mc A'' / X^\perp$ turning $(X',X)$ into a dual
   Banach algebra;
\item $X \subseteq \wap(\mc A')$;
\end{enumerate}
Furthermore, let $\mf A$ be a subalgebra of $\topcen{1}{\mc A''}
\cap \topcen{2}{\mc A''}$, and suppose that the natural map
$\mf A\rightarrow X'$ is surjective.  Then (1) and (2) hold,
and the algebra product given by (1) agrees with the product
induced by the map $\mf A\rightarrow X'$.
\end{proposition}
\begin{proof}
For $\Phi\in\mc A''$ and $\mu\in X$, suppose that $\Psi\in\mc A''$
is such that $\Phi+X^\perp = \Psi+X^\perp$, so that
\[ \ip{\Phi\cdot\mu}{a} = \ip{\Phi}{\mu\cdot a} = \ip{\Psi}{\mu\cdot a}
= \ip{\Psi\cdot\mu}{a} \qquad (a\in\mc A), \]
as $X$ is a submodule.  Hence there is a well-defined map
$(\mc A''/X^\perp) \times \mc A' \rightarrow \mc A'$ given by
$(\Phi + X^\perp) \cdot \mu = \Phi\cdot\mu$, and similarly with orders
reversed.  It is hence clear that $\aone$ gives a well-defined
product on $\mc A'' / X^\perp$ if and only if $\Phi\cdot\mu\in X$
for $\Phi\in\mc A''$ and $\mu\in X$.
If this holds, then $X$ is an $X'$-module if and only if,
for each $\mu\in X$ and $\Phi\in\mc A''$, there exists $\lambda\in X$
such that $\ip{\Phi\aone\Psi}{\mu} = \ip{\Psi}{\lambda}$ for
$\Psi\in\mc A''$.  In particular, we see that $\ip{\lambda}{a}
= \ip{\Phi\cdot a}{\mu} = \ip{\mu\cdot\Phi}{a}$ for $a\in\mc A$,
that is, $\lambda = \mu\cdot\Phi$, and so we conclude that
$\ip{\Phi\aone\Psi}{\mu} = \ip{\Phi\atwo\Psi}{\mu}$ for $\Phi,\Psi\in
\mc A'', \mu\in X$, which is equivalent to $X \subseteq \wap(\mc A')$.
Conversely, if (2) holds, then for $\Phi\in\mc A''$, $\mu\in X$ and
$\Psi\in X^\perp$, we have that $\ip{\Psi}{\Phi\cdot\mu} =
\ip{\Psi\atwo\Phi}{\mu} = \ip{\Phi}{\mu\cdot\Psi}=0$ as $\mu\cdot\Psi=0$,
which implies that $\Phi\cdot\mu\in X$.
Hence conditions (1) and (2) are equivalent.

Suppose that $\mf A\rightarrow X'$ is surjective,
so that for $\Phi,\Psi\in\mc A''$, let $a,b\in\mf A$ be
such that $a+X^\perp = \Phi+X^\perp$ and $b+X^\perp = \Psi+X^\perp$.
Then, for $\mu\in X$,
\begin{align*} \ip{\Phi\aone\Psi}{\mu} &= \ip{\Phi}{\Psi\cdot\mu}
= \ip{\Phi}{b\cdot\mu} = \ip{\Phi\aone b}{\mu} = \ip{\Phi\atwo b}{\mu}
= \ip{b}{\mu\cdot\Phi} \\ &= \ip{b}{\mu\cdot a} = \ip{a\atwo b}{\mu}
= \ip{a\aone b}{\mu}, \end{align*}
where we use the fact that $b\in\topcen{2}{\mc A''}$.  Thus the map
$\mf A\rightarrow X'$ gives a well-defined product on $X'$.
We also see that
\begin{align*} \ip{ab}{\mu} &= \ip{a}{b\cdot\mu} = \ip{a}{\Psi\cdot\mu}
= \ip{a\aone\Psi}{\mu} = \ip{a\atwo\Psi}{\mu} \\
&= \ip{\Psi}{\mu\cdot a} = \ip{\Psi}{\mu\cdot\Phi}
= \ip{\Phi\atwo\Psi}{\mu}, \end{align*}
as $a\in\topcen{1}{\mc A''}$.  Hence $\mu\in\wap(\mc A')$, and
(2) holds, as required.
\end{proof}

We note that in \cite[Proposition~4.9]{Kai}, Kaijser explores
similar ideas to the above proposition.  Furthermore, the
equivalence of (1) and (2) is established in \cite[Lemma~1.4]{LL2}
in the case of commutative Banach algebras.

\begin{corollary}
Let $\mc A$ be a Banach algebra.  Then $\wap(\mc A')'$ is a dual
Banach algebra.  Let $\mc A_* \subseteq \mc A'$ be a closed
submodule such that, if $\pi:\mc A'' \rightarrow \mc A''/\mc A_*^\perp
= \mc A_*'$ is the quotient map, then $\pi\circ\kappa_{\mc A}:
\mc A\rightarrow \mc A_*'$ is an isomorphism.  Then $\mc A_*'$
is a dual Banach algebra.
\end{corollary}

Conversely, suppose that $(\mc A,\mc A_*)$ is a dual Banach algebra.
Then it is a simple calculation (see \cite{Palmer}) that $\kappa_{\mc A_*}'
: \mc A'' \rightarrow \mc A$ is an algebra homomorphism for either
Arens product.  We may hence (and shall) make the following
equivalent definition:

\begin{definition}\label{other_defn}
Let $\mc A$ be a Banach algebra, let $\mc A_*$ be a closed
submodule of $\mc A'$, and let $\pi_{\mc A_*} : \mc A''
\rightarrow \mc A'' / \mc A_*^\perp = \mc A_*'$ be the quotient
map.  When $\pi_{\mc A_*} \circ \kappa_{\mc A} :\mc A\rightarrow
\mc A_*'$ is an isomorphism, we say that $\mc A$ is a
\emph{dual Banach algebra} with \emph{predual} $\mc A_*$.
\end{definition}

\begin{lemma}
Let $\mc A$ be a Banach algebra, and let $\mc B = \wap(\mc A')'$.
Then $\mc B$ is unital if and only if there exists $\Phi\in\mc A''$
with $\Phi\cdot\mu = \mu\cdot\Phi = \mu$ for each $\mu\in\wap(\mc A')$.
When $(\mc A,\mc A_*)$ is a dual Banach algebra, $\mc B$ is
unital if and only if $\mc A$ is unital.
\end{lemma}
\begin{proof}
Let $e_{\mc B}$ be the unit of $\mc B$, so that for some $\Phi\in\mc A''$,
we have that $\ip{\Phi}{\mu} = \ip{e_{\mc B}}{\mu}$ for $\mu\in\wap(\mc A')$.
Thus, for $\Psi\in\mc A''$,
\[ \ip{\Psi}{\mu} = \ip{\Psi\aone\Phi}{\mu}
= \ip{\Psi}{\Phi\cdot\mu} = \ip{\Phi\atwo\Psi}{\mu}
= \ip{\Psi}{\mu\cdot\Phi} \quad (\mu\in\wap(\mc A')), \]
as required.

Now suppose that $\mc A$ is a dual Banach algebra, and let $(a_\alpha)$
be a bounded net in $\mc A$ tending to $\Phi\in\mc A''$ in the
weak$^*$-topology.  Let $e\in\mc A$ be a weak$^*$-limit point of
$(a_\alpha)$.  Then, for $a\in\mc A$ and $\mu\in\mc A_*
\subseteq\wap(\mc A')$,
\[ \ip{ae}{\mu} = \ip{e}{\mu\cdot a}
= \lim_\alpha \ip{a_\alpha}{\mu\cdot a}
= \ip{\Phi}{\mu\cdot a} = \ip{a}{\mu}, \]
so that $ae=a$.  Similarly, $ea=a$, so that $e$ is a unit for $\mc A$.
\end{proof}

Thus looking at $\wap(\mc A')'$ is not useful for unitising a
dual Banach algebra; instead, $\wap(\mc A')'$ is a useful way
for embedding a Banach algebra in a dual Banach algebra.

\begin{lemma}
Let $\mc A$ and $\mc B$ be Banach algebras, and let $\pi:\mc A
\rightarrow\mc B$ be a homomorphism.  Then $\pi'(\wap(\mc B'))
\subseteq\wap(\mc A')$.
\end{lemma}
\begin{proof}
Let $\mu\in\wap(\mc B')$, and let $\lambda = \pi'(\mu)\in\mc A'$,
so that
\[ \ip{L_\lambda(a)}{b} = \ip{\lambda\cdot a}{b}
= \ip{\mu}{\pi(ab)} = \ip{\mu\cdot\pi(a)}{\pi(b)}
= \ip{\pi'L_\mu\pi(a)}{b}, \]
so that $L_\lambda = \pi'\circ L_\mu\circ \pi$ is weakly-compact,
as $L_\mu$ is weakly-compact.
\end{proof}

Weak$^*$-continuous representations of $\wap(\mc A')'$ are
closely related to continuous representations of $\mc A$, a fact
first noted by Runde.  The following is \cite[Theorem~4.10]{Runde2}.

\begin{proposition}\label{wap_reps}
Let $\mc A$ be a Banach algebra, let $(\mc B,\mc B_*)$ be
a dual Banach algebra, and let $\pi:\mc A\rightarrow\mc B$
be a homomorphism.  Then there is a unique weak$^*$-continuous
homomorphism $\hat\pi:\wap(\mc A')' \rightarrow\mc B$ such
that $\hat\pi\circ\kappa_{\mc A} = \pi$.  In particular, a
weak$^*$-continuous homomorphism $\theta:\wap(\mc A')'
\rightarrow\mc B$ is uniquely determined by its restriction to $\mc A$.
\end{proposition}

\begin{definition}\label{wap_alg_defn}
Let $\mc A$ be a Banach algebra.  We call $\wap(\mc A')'$ the
\emph{dual Banach algebra (DBA) enveloping algebra} of $\mc A$.
When the natural map of $\mc A$ into $\wap(\mc A')'$ is bounded
below, we say that $\mc A$ is a \emph{WAP algebra}.
\end{definition}

Then every dual Banach algebra, and every Arens regular Banach algebra,
is a WAP algebra.  We shall shortly see that group algebras are also
always WAP algebras, even in the non-discrete case, in which case they
are neither dual Banach algebras, nor Arens regular.  The above
proposition shows us that the weak$^*$-continuous theory of
the DBA enveloping algebra is determined by $\mc A$.
We hence see from this, and from later results, that the DBA
enveloping algebra plays much the same role as the enveloping
W$^*$-algebra of a C$^*$-algebra plays (see \cite[Chapter~III, Section~2]{tak}).

\section{Representations for WAP algebras}

We have already noted that work of Young shows that a WAP algebra
$\mc A$ admits a representation $\pi:\mc A\rightarrow\mc B(E)$
for some reflexive Banach space $E$, such that $\pi$ is bounded below.
In fact, Young effectively shows that $\mu\in\wap(\mc A')$ if
and only if there exists a reflexive Banach space $E$, a
representation $\pi:\mc A\rightarrow\mc B(E)$, $x\in E$ and
$\lambda\in E'$ with $\|x\| \|\lambda\| = \|\mu\|$ and such that
$\pi'\kappa_{E'\proten E}(\lambda\otimes x) = \mu$.  In particular,
we see that $\pi'\kappa_{E'\proten E}$ maps $E'\proten E$ onto
$\wap(\mc A')$.

In this section, we shall use some interpolation space theory
to prove an analogous result for dual Banach algebras which does
not seem to immediately follow from the results of Young (although
the method of proof is much the same).  We shall later use
interpolation space theory for other reasons, so it is useful to
define some concepts now.  Interpolation space arguments in this
area go back to \cite{DFJP}; we follow the text \cite{BB} for
results on interpolation spaces.

\begin{definition}\label{good_norm}
Let $(\mc A,\mc A_*)$ be a dual Banach algebra,
and let $\mu\in\mc A_*$.  Suppose that there exists a norm
$\|\cdot\|_\mu$ on $\mc A\cdot\mu = \{ a\cdot\mu : a\in\mc A\}$
such that the completion of $(\mc A\cdot\mu, \|\cdot\|_\mu)$,
denoted by $E_\mu$, is reflexive, and such that
\[ \|ab\cdot\mu\|_\mu \leq \|a\| \|b\cdot\mu\|_\mu,\qquad
\|a\cdot\mu\| \leq \|a\cdot\mu\|_\mu \leq \|a\|\|\mu\|
\qquad (a,b\in\mc A). \]
Let $\iota:E_\mu\rightarrow\mc A_*$ be the norm-decreasing
inclusion map, and suppose further that $\iota$ is injective.
Then we say that $\|\cdot\|_\mu$ is an \emph{admissible norm
for $\mu$}.
\end{definition}

\begin{example}
Let $(\mc A,\mc A_*)$ be a W$^*$-algebra, and let $\mu\in\mc A_*$
be a state.  Then it is simple to check that the usual GNS
construction for $\mu$ (see \cite[Chapter~I, Section~9]{tak})
induces an admissible norm on $\mc A\cdot\mu$.
\end{example}

\begin{example}
We note that $\|\cdot\|_\mu$ need not be unique (even in an isomorphic
sense).  For example, let $\mc A=l^2(\mathbb N)$ with pointwise multiplication,
and let $\mu\in \mc A'=l^2$ be such that the map $\mc A\rightarrow\mc A';
a\mapsto a\cdot\mu$ is injective and $\|\mu\|=1$ (for example,
$\mu=(2^{-n/2})_{n>0} \in l^2$).  Then define
\[ \|a\cdot\mu\|_{\mu,1} = \|a\cdot\mu\|, \quad \|a\cdot\mu\|_{\mu,2}
= \|a\| \qquad (a\in\mc A), \]
and let $E_{\mu,1}$ and $E_{\mu,2}$ be associated with $\|\cdot\|_{\mu,1}$
and $\|\cdot\|_{\mu,2}$ respectively.  Then $E_{\mu,1}$ is the closure
of $\mc A\cdot\mu$ in $l^2$, while $E_{\mu,2} = \mc A$, which are both
reflexive, as $l^2$ is reflexive.  We then check that
\begin{align*}
\|ab\cdot\mu\|_{\mu,1} &= \sup\{ |\ip{ab\cdot\mu}{c}| : \|c\|\leq 1\}
= \sup\{ |\ip{b\cdot\mu}{d}| : d=ca, \|c\|\leq 1 \} \\
&\leq \sup\{ |\ip{b\cdot\mu}{d}| : \|d\| \leq \|a\| \}
= \|a\| \|b\cdot\mu\| = \|a\| \|b\cdot\mu\|_{\mu,1},
\end{align*}
while clearly $\|ab\cdot\mu\|_{\mu,2} \leq \|a\| \|b\cdot\mu\|_{\mu,2}$.
Hence both $\|\cdot\|_{\mu,1}$ and $\|\cdot\|_{\mu,2}$ are admissible,
but clearly they are not equivalent norms.
\end{example}

\begin{lemma}
Let $(\mc A,\mc A_*)$ be a dual Banach algebra, let $\mu\in\mc A_*$ have
an admissible norm, and let $E_\mu$ be a space as defined above using some
admissible norm for $\mu$.  Then $\iota':\mc A\rightarrow E_\mu'$ has
dense range, the module action of $\mc A$ on $E_\mu$ induces a
weak$^*$-continuous representation $\mc A\rightarrow\mc B(E_\mu)$, and
there exists $x\in E_\mu$ and $\lambda\in E_\mu'$ such that
$\|x\|\|\lambda\|=\|\mu\|$ and $\ip{\lambda}{a\cdot x} = \ip{a}{\mu}$
for $a\in\mc A$.
\end{lemma}
\begin{proof}
We may suppose that $\|\mu\|=1$.  By the condition on $\|\cdot\|_\mu$,
we see that $\iota$ is norm-decreasing, and the module action
inherited from $\mc A_*$ induces a Banach left $\mc A$-module
action on $E_\mu$.  Furthermore, $\iota'$ has dense range if
and only if $\iota'':E_\mu''\rightarrow\mc A'$ is injective, which,
as $E_\mu$ is reflexive, is in turn equivalent to $\iota$ being
injective.

We define $\psi_\mu : E_\mu' \proten E_\mu \rightarrow \mc A_*$ by
\[ \psi_\mu\big( \iota'(a) \otimes b\cdot\mu \big) = b\cdot\mu\cdot a
\qquad (a\in\mc A, b\cdot\mu\in E_\mu). \]
Assuming this is bounded, $\psi_\mu$ extends by linearity and continuity
to $E_\mu' \proten E_\mu$.  Indeed, we have
\begin{align*}
\|b\cdot\mu\cdot a\| &= \sup\{ |\ip{acb}{\mu}| : \|c\|\leq 1 \}
= \sup\{ |\ip{a}{cb\cdot\mu}| : \|c\|\leq 1 \} \\
&\leq \sup\{ |\ip{a}{d\cdot\mu}| : \|d\cdot\mu\|_{E_\mu}\leq \|b\cdot\mu\|_{E_\mu} \}
= \| \iota'(a) \|_{E_\mu'} \|b\cdot\mu\|_{E_\mu},
\end{align*}
as $\|cb\cdot\mu\|_{E_\mu}\leq \|c\|\|b\cdot\mu\|_{E_\mu}$.
Thus $\psi_\mu$ is norm-decreasing.  Then let $\theta_\mu =
\psi_\mu' : \mc A\rightarrow \mc B(E_\mu)$, so that for
$a,b,c\in\mc A$,
\[ \ip{\iota'(b)}{\theta_\mu(a)(c\cdot\mu)}
= \ip{a}{\psi_\mu( \iota'(b) \otimes c\cdot\mu )}
= \ip{b}{ac\cdot\mu} = \ip{\iota'(b)}{a\cdot(c\cdot\mu)}, \]
so that $\theta_\mu$ agrees with the left-module action of $\mc A$
on $E_\mu$, as required.

Finally, we see that $\psi_\mu(\iota'(e_{\mc A}) \otimes e_{\mc A}\cdot\mu)
= \mu$, where
\begin{align*}
\|\iota'(e_{\mc A})\|_{E_\mu'} &= \sup\{ |\ip{e_{\mc A}}{a\cdot\mu}| :
   \| a\cdot\mu \|_{E_\mu}\leq 1 \} \\
&\leq \sup\{ |\ip{e_{\mc A}}{a\cdot\mu}| : \| a\cdot\mu \|\leq 1 \}
\leq \|e_{\mc A}\| = 1,
\end{align*}
and $\|e_{\mc A}\cdot\mu\|_{E_\mu} \leq \|e_{\mc A}\| \|\mu\| = 1$.
\end{proof}

\begin{theorem}\label{rep_thm}
Let $(\mc A,\mc A_*)$ be a dual Banach algebra such that each norm-one
member of $\mc A_*$ has an admissible norm.  Then $\mc A$ is isometric,
via a weak$^*$-weak$^*$-continuous map, to a weak$^*$-closed subalgebra
of $\mc B(E)$ for some reflexive Banach space $E$.
\end{theorem}
\begin{proof}
We may suppose that $\mc A$ is unital, as otherwise, we may
work with $(\mc A,\mc A_*)^\sharp$, and then restrict the resulting
representation to $\mc A$.
Let $X = \{\mu\in\mc A_* : \|\mu\|=1 \}$, and let
$E = l^2\big(\bigoplus_{\mu\in X} E_\mu\big)$,
so that $E$ is a reflexive Banach space.  Define $\psi:E'\proten
E \rightarrow \mc A_*$ by
\[ \psi\big( (\lambda_\mu) \otimes (x_\mu) \big)
= \sum_{\mu\in X}\psi_\mu(\lambda_\mu\otimes x_\mu)
\qquad ((\lambda_\mu)\in E', (x_\mu)\in E). \]
This is norm-decreasing, as
\[ \Big\| \sum_{\mu\in X} \psi_\mu(\lambda_\mu\otimes x_\mu) \Big\|
\leq \sum_{\mu\in X} \|\lambda_\mu\| \|x_\mu\|
\leq \Big( \sum_{\mu\in X} \|\lambda_\mu\|^2 \Big)^{1/2}
\Big( \sum_{\mu\in X} \|x_\mu\|^2 \Big)^{1/2}. \]
Then let $\theta = \psi' : \mc A \rightarrow \mc B(E)$, so
that $\theta$ is weak$^*$-continuous, and for $a,b\in\mc A$,
$(x_\mu)\in E$ and $(\lambda_\mu)\in E'$, we have
\begin{align*}
\ip{(\lambda_\mu)}{\theta(b)(x_\mu)}
&= \sum_{\mu\in X} \ip{b}{\psi_\mu(\lambda_\mu\otimes x_\mu)}
= \sum_{\mu\in X} \ip{\lambda_\mu}{\theta_\mu(b)(x_\mu)},
\end{align*}
so that $\theta(b)(x_\mu) = \big( \theta_\mu(b)(x_\mu) \big)$,
and hence $\theta$ is a homomorphism, as each $\theta_\mu$ is
a homomorphism.

It is a standard result that $\theta=\psi'$ has a weak$^*$-closed
image if and only if $\theta$ has a closed image, which is if and
only if $\theta$ is bounded below.  For $a\in\mc A$ and
$\epsilon>0$, there exists $\mu\in X$ such that $|\ip{a}{\mu}| >
(1-\epsilon)\|a\|$.  Then $\mu = \psi(\lambda_\mu\otimes x_\mu)$
for some $\lambda_\mu\in E_\mu'$ and $x_\mu\in E_\mu$ with
$\|\lambda_\mu\| \|x_\mu\|=1$.  Thus
\[ \|\theta(a)\| \geq \|\theta_\mu(a)\| \geq
|\ip{a}{\psi(\lambda_\mu\otimes x_\mu)}| > (1-\epsilon)\|a\|. \]
As $\epsilon>0$ was arbitrary, we see that $\theta$ is an isometry
onto its range.
\end{proof}

Of course, we haven't shown that any dual Banach algebra
(other than a W$^*$-algebra) admits such a representation.  We
now remedy this situation by using some interpolation space theory.

\begin{theorem}\label{int_thm}
Let $\mc A$ be a unital dual Banach algebra with predual $\mc A_*$.
Then each norm-one member of $\mc A_*$ has an admissible norm.
\end{theorem}
\begin{proof}
Let $\mu\in\mc A_*$ be such that $\|\mu\|=1$, and for
$n\in\mathbb N$ define a new norm on $\mc A'$ by
\[ \|\lambda\|_n = \inf\{ 2^{-n/2} \|b\| + 2^{n/2}\|\lambda - b\cdot\mu\|
: b\in\mc A \} \qquad (\lambda\in\mc A'). \]
Then, for $a,b\in\mc A$ and $\lambda\in\mc A'$, we have that
$\|\lambda\|_n \leq 2^{n/2}\|\lambda\|$, $\|a\cdot\mu\|_n \leq
2^{-n/2} \|a\|$, and $\|\lambda\| \leq
\|\lambda - b\cdot\mu\| + \|b\cdot\mu\| \leq 2^n \|\lambda-b\cdot\mu\|
+ \|b\|$ so that $2^{-n/2}\|\lambda\| \leq \|\lambda\|_n$.
We then define
\[ E_\mu = \Big\{ \lambda\in\mc A' : \|\lambda\|_\mu :=
\Big(\sum_{n=1}^\infty \|\lambda\|_n^2 \Big)^{1/2} < \infty \Big\}. \]
If $\lambda\in E_\mu$, then there exists some sequence $(b_n)$ in $\mc A$
such that $2^n\|\lambda-b_n\cdot\mu\|^2 \rightarrow 0$ as $n\rightarrow\infty$.
In particular, $\mc A\cdot\mu$ is dense in $E_\mu$.
Thus, for $a\in\mc A$, we have that
$\|a\cdot\mu\| \leq \|a\cdot\mu\|_\mu \leq \|a\|$.
We can also easily check that $\|a\cdot\lambda\|_\mu \leq
\|a\| \|\lambda\|_\mu$ for $a\in\mc A$ and $\lambda\in E_\mu$.
Hence we need only show that $E_\mu$ is reflexive to verify the
conditions of Definition~\ref{good_norm}.

Recall the definition of $R_\mu:\mc A\rightarrow\mc A'$, that is,
$R_\mu(a)=a\cdot\mu$.  Then $R_\mu$ maps into $\mc A_*$ and is
weakly-compact, as $\mc A*\subseteq\wap(\mc A')$.  It follows from the work
in \cite{DFJP} (see \cite[Section~1.7.8]{PalBook} for a sketch) that
$E_\mu$ is reflexive, as the map $R_\mu$ is weakly-compact.
\end{proof}

Notice that the above proof will work for any $\mu\in\wap(\mc A')$,
which re-creates Young's result.

\begin{remark}
The above construction of $E_\mu$ is actually a Lions-Peetre
interpolation space.  Let $\mc A\cdot\mu$ be the subspace of
$\mc A_*$ spanned by $\{ a\cdot\mu : a\in\mc A\}$ together
with the norm $\|a\cdot\mu\|_{\mc A\cdot\mu} =
\inf\{ \|b\| : b\cdot\mu = a\cdot\mu \}$.  Then we see that
$R_\mu:\mc A\rightarrow \mc A\cdot\mu$ is norm-decreasing, and
the induced map $\mc A / \ker R_\mu \rightarrow \mc A\cdot\mu$
is an isometry, showing that $\mc A\cdot\mu$ is a Banach space.

Following \cite{BB}, we let $\mc S(\mc A\cdot\mu, \mc A_*)$ be
the space $\mc A_*$ together with the norm
\[ \|\lambda\|_{\mc S} = \inf\big\{ \|a\cdot\mu\|_{\mc A\cdot\mu}
+ \|\phi\| : \lambda = a\cdot\mu + \phi \big\}
\qquad (\lambda\in\mc A_*). \]
Then $\mc S$ is a Banach space.  Let $1\leq p<\infty$, $\xi_0<0$ and
$\xi_1>0$, and let $s_2^+(p; \xi_0,\mc A\cdot\mu; \xi_1,\mc A_*)$ be the
subspace of $\mc S$ such that
\begin{align*}
\|\lambda\|_{s_p^+} &= \inf\Big\{ \max\Big(
\Big( \sum_{n=1}^\infty \|e^{\xi_0 n}a_n\cdot\mu\|_{\mc A\cdot\mu}^p \Big)^{1/p},
\Big( \sum_{n=1}^\infty \|e^{\xi_1 n}\phi_n\|^p \Big)^{1/p} \Big) \\
&\hspace{15ex}: \lambda = a_n\cdot\mu + \phi_n\, (n>0) \Big\} < \infty.
\end{align*}
In comparison, we see that $E_\mu$ is the subspace of $\mc S$ such that
\[ \|\lambda\|_\mu = \inf\Big\{ \Big( \sum_{n=1}^\infty
\big( 2^{-n/2}\|a_n\cdot\mu\|_{\mc A\cdot\mu} +
2^{n/2}\|\phi_n\| \big)^2 \Big)^{1/2} :
\lambda = a_n\cdot\mu + \phi_n\,(n>0) \Big\}. \]
We thus see that $s_2^+(2; -\log \sqrt 2,\mc A\cdot\mu; \log\sqrt 2,
\mc A_*)$ is isomorphic to $E_\mu$.

It follows from \cite[Section~1.5, Proposition~1]{BB} that
$E_\mu$ is a member of the isomorphic class
$(\mc A\cdot\mu, \mc A_*)_{1/2,2}$, and hence
\cite[Section~2.3, Proposition~1]{BB} tells us that $E_\mu$
is reflexive if and only if the inclusion $\mc A\cdot\mu
\rightarrow\mc A_*$ is weakly-compact, which is if and only if
$R_\mu$ is weakly-compact, as before.
\end{remark}

Putting all these results together, we obtain the following.

\begin{corollary}
Let $(\mc A,\mc A_*)$ be a dual Banach algebra.  Then $\mc A$
admits an isometric, weak$^*$-weak$^*$-continuous representation
on some reflexive Banach space.
\end{corollary}

\begin{proposition}
Let $(\mc A,\mc A_*)$ be a W$^*$-algebra (that is, $\mc A$
is a C$^*$-algebra), and form $E$ as above using only the
$\mu\in\mc A_*$ which are states.  Then $E$ is isomorphic
to a Hilbert space, and our representation agrees with the usual
universal representation for a W$^*$-algebra.
\end{proposition}
\begin{proof}
This follows from work in \cite{DawLau}.
\end{proof}

\section{Unique preduals}

It is a standard result in the theory of W$^*$-algebras that
a W$^*$-algebra has a unique predual, up to isometric
classification.  That is, if $(\mc A,\mc A_*)$ is a W$^*$-algebra,
$E$ is a Banach space, and $\theta:\mc A\rightarrow E'$ is an
isometric isomorphism, then $\theta$ is automatically
weak$^*$-continuous.  This follows as we can use $\theta$ to
induce a C$^*$-algebra structure on $E'$, showing that $E$ is
also a predual for $\mc A$.

The theory of isometric preduals in Banach spaces has
attracted some attention (see the survey \cite{Gode}).  However,
here we are interested in the \emph{isomorphic} not \emph{isometric}
theory.

\begin{theorem}
Let $(\mc A,\mc A_*)$ be a commutative W$^*$-algebra, let
$(\mc B,\mc B_*)$ be a dual Banach algebra, and let
$\theta:\mc A\rightarrow \mc B$ be a Banach algebra isomorphism.
Then $\theta$ is automatically weak$^*$-continuous.
\end{theorem}
\begin{proof}
By \cite[Theorem~1.18]{tak}, $\mc A$ can be identified with $C(\Omega)$
where $\Omega$ is the character space of $\mc A$, which is a
\emph{hyperstonian} space (see \cite[Section~8]{Bade} for
further details).  In particular, $\Omega$ is \emph{stonian}
in that the closure of any open set is open.  Then $\mc A'$ is
$M(\Omega)$, the space of regular Borel measures on $\Omega$.
We say that a positive measure $\mu\in M(\Omega)$ is \emph{normal}
if $\mu(A)=0$ whenever $A$ is \emph{nowhere dense}, that is, the
closure of $A$ has empty interior.  A general measure is \emph{normal}
when its absolute value is normal.
By \cite[Theorem~8.2]{Bade} this definition agrees with the
usual one for W$^*$-algebras (see also \cite[Proposition~1.11]{tak}).
Then the collection of normal measures forms a closed subspace of
$M(\Omega)$ which is equal to $\kappa_{\mc A_*}(\mc A_*)$.

By reversing the argument which lead to Definition~\ref{other_defn},
we need to show that if $E\subseteq\mc A'$ is a closed submodule such
that the natural map $\iota_E:\mc A\rightarrow E' = \mc A''/E^\perp$
is an isomorphism, then $E = \kappa_{\mc A_*}(\mc A_*)$.  The
equivalence of this statement to the statement involving $\mc B$
comes from setting $E = \theta'\kappa_{\mc B_*}(B_*)$.

Let $E$ be as stated, and choose $\lambda\in E$.  We will show
that $\lambda$ is normal, which will complete the proof, as then
$E \subseteq \mc A_*$, and so necessarily $E = \mc A_*$.  Let $A$
be a closed subset of $\Omega$ with empty interior, so we aim to
show that $|\lambda|(A)=0$.  Consider the family $\mc C$ of closed
and open subsets of $\Omega$ which contain $A$, partially ordered
by reverse inclusion.  For $B\in\mc C$, let $\chi_B$ be the indicator
function of $B$, so that $\chi_B \in C(\Omega)$.  As $E$ is a predual,
there exists a unique $f\in C(\Omega)$ such that
\[ \ip{\mu}{f} = \lim_{B\in\mc C} \ip{\mu}{\chi_B}
\qquad (\mu\in E). \]
For $C\in\mc C$, notice that
\[ \ip{\mu}{\chi_C f} = \ip{\mu\cdot\chi_C}{f}
= \lim_{B\in\mc C} \ip{\mu\cdot\chi_C}{\chi_B}
= \lim_{B\in\mc C} \ip{\mu}{\chi_B} = \ip{\mu}{f}
\quad (\mu\in E), \]
so that $\chi_C f = f$, and hence, for $x\in\Omega\setminus C$,
$f(x) = f(x)\chi_C(x) = 0$.  We hence see that $f$ vanishes off the set
$A_0 := \bigcap_{B\in\mc C} B$.

We claim that $A_0 = A$, which follows from some simple topology.
Indeed, clearly $A\subseteq A_0$, and suppose towards a contradiction
that there exists $x\in A_0 \setminus A$.  Thus, for each open
$B\subseteq\Omega$ with $A\subseteq\overline B$, we have that $x\in\overline B$.
As $\Omega$ is Hausdorff, for each $a\in A$ there exist disjoint
open sets $U_a$ and $V_a$ with $a\in U_a$ and $x\in V_a$.  As $A$
is compact, there exists $a_1,\ldots,a_n$ in $A$ such that
$A\subseteq U:= U_{a_1}\cup\cdots\cup U_{a_n}$.  Clearly $U$ is
disjoint from the open set $V:= V_{a_1}\cap\cdots\cap V_{a_n}$,
so that as $x\in V$, we see that $x\not\in\overline{U}$, a contradiction.

Consequently $f$ is supported on $A$, and as $A$ has empty interior and
$f$ is continuous, we must have that $f=0$.
Let $\lambda = \lambda_r +i\lambda_i$ where
$\lambda_r$ and $\lambda_i$ are real measures.  Then, as
$\lim_{B\in\mc C}\ip{\lambda}{\chi_B}=0$, we see that
\[ \lim_{B\in\mc C} \ip{\lambda_r}{\chi_B} =
\lim_{B\in\mc C} \ip{\lambda_i}{\chi_B} = 0, \]
as $\chi_B$ is real-valued.  There exists a \emph{Hahn decomposition}
see \cite[Section~29]{Halmos} for $\lambda_r$, that is, measurable sets
$E_+$ and $E_-$ such that $\Omega = E_+ \cup E_-$ and with
\[ \lambda_+(E) = \lambda_r(E\cap E_+),
\quad \lambda_-(E) = - \lambda_r(E\cap E_-) \qquad (E\subseteq\Omega), \]
defining two positive measures $\lambda_+$ and $\lambda_-$ with
$|\lambda_r| = \lambda_+ + \lambda_-$ and $\lambda_r = \lambda_+ -
\lambda_-$.  As $|\lambda_r|$ is regular, for each $\epsilon>0$,
there exists an open set $U$ and a closed set $K$ such that
$K \subseteq E_+ \subseteq U$ with $|\lambda_r|(U\setminus K)< \epsilon$.
As $\Omega$ is Stonian, we can find an open and closed set $V$ such
that $K \subseteq V \subseteq U$ (this follows by a similar argument
to that employed above to show that $A_0=A$).  Then
\[ |\lambda_r|(E_+ \setminus V) + |\lambda_r|(V\setminus E_+)
\leq |\lambda_r|(U\setminus K) < \epsilon. \]
We hence see that for $B\subseteq\Omega$,
\begin{align*}
\lambda_+ (B) &= \lambda(B\cap E_+)
= \big| \lambda_r(B\cap V) + \lambda_r(B\cap(E_+\setminus V))
- \lambda_r(B\cap(V\setminus E_+)) \big| \\
&\leq \epsilon + |\lambda_r(B\cap V)|.
\end{align*}
Consequently, 
\[ 0 \leq \lim_{B\in\mc C} \ip{\lambda_+}{\chi_B}
\leq \epsilon + \lim_{B\in\mc C} \ip{\lambda_r}{\chi_V\chi_B}
= \epsilon + \lim_{B\in\mc C} \ip{\lambda_r\cdot \chi_V}{\chi_B}
= \epsilon, \]
so as $\epsilon>0$ was arbitrary, $\lim_{B\in\mc C}
\ip{\lambda_+}{\chi_B}=0$.  We then see that
\[ 0\leq \lambda_+(A) \leq \inf_{B\in\mc C} \lambda_+(B)
= 0. \]
A similar argument shows that $\lambda_-(A)=0$, so that
$|\lambda_r|(A)=0$.  Similarly, $|\lambda_i|(A)=0$, so that
$|\lambda|(A)=0$.  We thus see that $\lambda$ is normal, as required.
\end{proof}

We note that we cannot drop the assumption that $\theta$ is
an algebra homomorphism.  This follows as, for example,
Pe\l czy\'nski showed in \cite{Pel} that $l^\infty$ and
$L^\infty[0,1]$ are isomorphic (but not isometric, and
not isomorphic as Banach algebras in their natural
products).  However, of course, $l^1$ and $L^1[0,1]$
are not isomorphic, so no isomorphism $l^\infty
\rightarrow L^\infty[0,1]$ can be weak$^*$-continuous.

Continuing the theme of unique preduals, we have the following,
which is \cite[Proposition~5.10]{GS}.

\begin{proposition}
Let $E$ be a reflexive Banach space.  Then the predual $E'\proten E$
is isometrically unique for $\mc B(E)$, meaning that when
$\phi:F'\rightarrow\mc B(E)$ is an isometric isomorphism for some
Banach space $F$, there exists an isometry $\psi:E'\proten E
\rightarrow F$ such that $\psi'=\phi$.
\end{proposition}

\begin{definition}
Let $E$ be a Banach space, and suppose that for each compact
set $K\subseteq E$ and each $\epsilon>0$, there exists a finite-rank
operator $T\in\mc F(E)$ such that $\|T(x)-x\|<\epsilon$ for
each $x\in K$.  Then we say that $E$ has the \emph{approximation
property}.  When we can choose $T$ to be uniformly bounded, $E$ has
the \emph{bounded approximation property}, and when we can choose $T$
be a contraction, $E$ has the \emph{metric approximation property}.
\end{definition}

See \cite[Chapter~4]{Ryan} or \cite[Chapter~VIII, Section~3]{DU}
for further details.  When $E$ is a reflexive Banach space, the
approximation property implies the metric approximation property
(see \cite[Corollary~5.51]{Ryan}).  The approximation property is
equivalent to the natural map
$E'\proten E\rightarrow\mc B(E)$ being injective.

\begin{theorem}
Let $E$ be a reflexive Banach space with the approximation property.
Whenever $(\mc A,\mc A_*)$ is a dual Banach algebra and
$\theta:\mc B(E)\rightarrow\mc A$ is a Banach algebra isomorphism,
$\theta$ is weak$^*$-continuous.
\end{theorem}
\begin{proof}
Again, we need to show that when $X\subseteq\mc B(E)'$ is a
predual, we have that $X = \kappa_{E'\proten E}(E'\proten E)$.
As $E$ is reflexive and has the approximation property, $\mc A(E)'
= E'\proten E$ and so $\mc A(E)'' = \mc B(E)$.  We may check that
$\mc A(E)$ is Arens regular, and that the Arens products on
$\mc A(E)''$ agree with the usual product on $\mc B(E)$.  By
standard results, as $\mc B(E)$ is unital, there exists a bounded
approximate identity $(e_\alpha)$ for $\mc A(E)$ (see
\cite[Proposition~2.9.16]{Dales}
or \cite[Section~1.7.13]{PalBook} for further details).
Let $\iota:\mc A(E)\rightarrow\mc B(E)$ be the inclusion map,
so that actually $\iota$ agrees with the map $\kappa_{\mc A(E)}$.
Then $\iota':\mc B(E)'\rightarrow E'\proten E$ satisfies
$\iota'\kappa_{E'\proten E} = I_{E'\proten E}$. As $X$ is a predual,
there exists a unique $P\in\mc B(E)$ such that
\[ \ip{\mu}{P} = \lim_\alpha \ip{\mu}{\iota(e_\alpha)}
= \lim_\alpha \ip{\iota'(\mu)}{e_\alpha}
= \ip{I_E}{\iota'(\mu)} \qquad (\mu\in X). \]
Then, for $T\in\mc B(E)$ and $\mu\in X$,
\begin{align*} \ip{\mu}{PT} &= \ip{T\cdot\mu}{P} =
\lim_\alpha \ip{T\cdot\mu}{\iota(e_\alpha)}
= \lim_\alpha \ip{\mu}{\iota(e_\alpha T)}
= \lim_\alpha \ip{\iota'(\mu)}{e_\alpha T} \\
&= \lim_\alpha \ip{T\cdot\iota'(\mu)}{e_\alpha}
= \ip{I_E}{T\cdot\iota'(\mu)} = \ip{T}{\iota'(\mu)} \end{align*}
and so, for each $\alpha$,
\[ \ip{\mu}{P\iota(e_\alpha)} = \ip{\iota(e_\alpha)}{\iota'(\mu)}
= \ip{\mu}{\iota(e_\alpha)}  \qquad (\mu\in X), \]
so that $P\iota(e_\alpha) = \iota(e_\alpha)$.  Similarly,
$\iota(e_\alpha)P = \iota(e_\alpha)$ for every $\alpha$.  Hence
\[ \ip{\mu}{P^2} = \lim_\alpha \ip{\mu\cdot P}{\iota(e_\alpha)}
= \lim_\alpha \ip{\mu}{P\iota(e_\alpha)}
= \lim_\alpha \ip{\mu}{\iota(e_\alpha)} = \ip{\mu}{P} \]
for $\mu\in X$, so that $P$ is a projection.  As $( P\iota(e_\alpha) )
= ( \iota(e_\alpha) )$ is a bounded approximate identity
for $\mc A(E)\subseteq\mc B(E)$, we must have that the image of
$P$ is the whole of $E$, that is, $P$ is the identity.  Thus
\[ \ip{\mu}{T} = \ip{\mu}{PT} = \ip{T}{\iota'(\mu)}
= \ip{\kappa_{E'\proten E}\iota'(\mu)}{T}
\qquad (T\in\mc B(E), \mu\in X), \]
and hence $\kappa_{E'\proten E}\iota'$ is the identity
on $X$.  This implies that $X \subseteq \kappa_{E'\proten E}
(E'\proten E)$ so that $X = \kappa_{E'\proten E}(E'\proten E)$,
as required.
\end{proof}

It would be nice to remove the condition on $E$ having the
approximation property, but this is utterly integral to the current
proof.  The general question of which dual Banach algebras have
a unique predual seems very interesting.  We have looked
at, but have been unable to answer, the question of whether $\ell^1(\mathbb Z)$
has a unique predual.  This is equivalent to the very concrete
question: let $X \subseteq \ell^\infty(\mathbb Z)$ be a shift-invariant
subspace such that $X'$ is naturally identified with $\ell^1(\mathbb Z)$.
Is $X = c_0(\mathbb Z)$?  Of course, as a Banach space, $\ell^1(\mathbb Z)$
has plenty of preduals (see \cite{Gode}) but these do not appear to
be preduals which make the product on $\ell^1(\mathbb Z)$ weak$^*$-continuous.

\section{Dual Banach $*$-algebras}

We start by studying dual Banach algebras $(\mc A,\mc A_*)$ which admit
an involution, which for us shall be a continuous, conjugate-linear
map $*:\mc A\rightarrow\mc A$ such that $(ab)^* = b^* a^*$.
Recall that we may define an involution $*$ on $\mc A'$ by
\[ \ip{\mu^*}{a} = \overline{ \ip{\mu}{a^*} }
\qquad (\mu\in\mc A', a\in\mc A). \]
Then the involution is weak$^*$-continuous if and only if
$\mc A_*$ forms a self-adjoint subspace of $\mc A'$.

We shall now sketch some results on Banach spaces which admit
a sesquilinear form which is not necessarily positive (such spaces
are hence generalisations of Hilbert spaces).  These are studied by
Laustsen and the author in \cite{DawLau}.  Let $E$ be a Banach
space and let $[\cdot,\cdot]$ be a sesquilinear form on $E$ which is
\emph{bounded} in the sense that for some $C>0$, $|[x,y]| \leq C
\|x\| \|y\|$ for $x,y\in E$.  There hence exists a bounded,
conjugate-linear map $J:E\rightarrow E'$ such that $[x,y]
= \ip{J(y)}{x}$ for each $x,y\in E$.  Suppose that $J$ is a
homeomorphism (which forces $E$ to be reflexive).  Then there is
an involution on $\mc B(E)$ given equivalently by
\[ T^* = J^{-1} \circ T' \circ J, \quad\text{or}\quad [T^*(x),y] = [x,T(y)]
\qquad (x,y\in E, T\in\mc B(E)). \]
It is shown in \cite{DawLau} that every bounded involution on
$\mc B(E)$ arises in this way.

Now let $(\mc A,\mc A_*)$ be a dual Banach algebra with a
continuous involution, let $\mu\in\mc A_*$ be self-adjoint and
of norm one, and consider the space $E_\mu$ formed in
Theorem~\ref{int_thm}.  As this space is isomorphic to a Lions-Peetre
interpolation space, it is isomorphic to the spaces constructed in
\cite{DawLau}.  In particular, we may define a bounded form on $E_\mu$ by
\[ [a\cdot\mu, b\cdot\mu] = \ip{b^*a}{\mu} \qquad (a,b\in\mc A) \]
such that the induced map $J:E\rightarrow E'$ is a homeomorphism.
It hence follows that the representation $\mc A\rightarrow\mc B(E)$
is actually a $*$-homomorphism.

The space $(E,[\cdot,\cdot])$ may certainly fail to be a Hilbert space.
However, it is shown that in the special case when $\mu$ is a positive
linear functional, then $E$ is at least isomorphic to the Hilbert
space generated by the GNS representation for $\mu$.

\begin{proposition}\label{star_has_rep}
Let $(\mc A,\mc A_*)$ be a dual Banach algebra with a continuous
involution.  Then the following are equivalent:
\begin{enumerate}
\item the involution is weak$^*$-continuous;
\item $\mc A$ is weak$^*$-continuously $*$-isomorphic to a
closed subspace of $\mc B(E)$ for some reflexive Banach space $E$
such that $\mc B(E)$ admits an involution.
\end{enumerate}
\end{proposition}
\begin{proof}
To show that (1) implies (2), by the preceding discussion,
the only thing to check is that the subset of $\mc A_*$ of norm-one
self-adjoint functionals norms $\mc A$.  However, as the involution
is weak$^*$-continuous, $\mc A_*$ is itself self-adjoint, and
hence every $\mu\in\mc A_*$ is of the form $\mu = \mu_r + i\mu_i$
for self-adjoint $\mu_r,\mu_i\in\mc A_*$.  It hence follows that
the quantity $\sup\{ |\ip{a}{\mu}| : \mu\in\mc A_*, \mu^*=\mu, \|\mu\|=1 \}$
is at least equivalent to $\|a\|$, as required.

Now suppose that (2) holds.  We can hence identify $\mc A$ with
its image in $\mc B(E)$, so that $\mc A_*$ becomes identified
with $E' \proten E / {^\perp}\mc A_*$.
Let $\tau\in\mc A_*$, so that as $J$ is a homeomorphism, $\tau$ has
a representation of the form
\[ \tau = \sum_{n=1}^\infty J(x_n) \otimes y_n + {^\perp}\mc A_*, \]
with $\sum_{n=1}^\infty \|x_n\| \|y_n\|<\infty$.  Let
\[ \sigma = \sum_{n=1}^\infty J(y_n) \otimes x_n + {^\perp}\mc A_*
\in \mc A_*. \]
Then, for $T\in\mc A \subseteq\mc B(E)$, we see that
\[ \ip{T}{\tau} = \sum_{n=1}^\infty \ip{J(x_n)}{T(y_n)}
= \sum_{n=1}^\infty [ T(y_n), x_n ]
= \sum_{n=1}^\infty \overline{ [ T^*(x_n), y_n ] }
= \overline{ \ip{T^*}{\sigma} }, \]
so that $\ip{T}{\sigma^*} = \overline{ \ip{T^*}{\sigma} }
= \ip{T}{\tau}$.  Hence $\sigma^* = \tau$, so that $\tau^* = \sigma$,
and hence $\mc A_*$ is self-adjoint, as required.
\end{proof}

It would be interesting to know if there exists a dual Banach
algebra $(\mc A,\mc A_*)$ which admits a continuous involution
which is not weak$^*$-continuous.

\begin{lemma}\label{eq_inv_norm}
Let $(\mc A,\mc A_*)$ be a dual Banach algebra, and suppose
that $\mc A$ admits a weak$^*$-continuous involution.
Then there is an equivalent norm $\|\cdot\|_0$ on
$\mc A_*$ such that the involution becomes isometric
on $(\mc A,\|\cdot\|_0^*)$, where $\|\cdot\|_0^*$ is the
dual norm to $\|\cdot\|_0$.
\end{lemma}
\begin{proof}
We would usually define a new norm on $\mc A$ by
$|||a||| = \max(\|a\|, \|a^*\|)$; we show here how to
dualise this idea.  For $\mu\in A_*$, we define
\[ \|\mu\|_0 = \inf \{ \|\lambda^*\| + \|\mu-\lambda\|
: \lambda\in\mc A_* \}, \]
where, as $\mc A_*$ is self-adjoint, $\lambda^* \in\mc A_*$.
Clearly we have that $\|\mu\| \leq \|\mu\|_0$ for $\mu\in\mc A_*$.
As the involution is continuous, there exists $M\geq 1$ such that
$\|a^*\|\leq M\|a\|$ for each $a\in\mc A$.  Then $\| \lambda^* \|
\leq M\|\lambda\|$ for each $\lambda\in\mc A'$, so that also
$\| \lambda \| \leq M\|\lambda^*\|$.  Thus
\[ \|\mu\| \leq \|\lambda\| + \|\mu-\lambda\|
\leq M\|\lambda^*\| + \|\mu-\lambda\|
\leq M\|\lambda^*\| + M\|\mu-\lambda\|, \qquad (\lambda\in\mc A_*) \]
so that $\| \mu \| \leq M \| \mu \|_0$.
Hence $\|\cdot\|_0$ is equivalent to $\|\cdot\|$ on $\mc A_*$.

Then, for $a\in\mc A$, we have that
\begin{align*}
\|a\|_0^* &= \sup\{ |\ip{a}{\mu+\lambda}| : \mu,\lambda\in\mc A_*,
   \|\lambda^*\|+\|\mu\|\leq 1 \} \\
&= \sup\{ |\ip{a}{\mu}| + |\ip{a}{\lambda}| : \mu,\lambda\in\mc A_*,
   \|\lambda^*\|+\|\mu\|\leq 1 \} \\
&= \max\big\{ \sup\{ |\ip{a}{\mu}| : \|\mu\|\leq 1 \},
   \sup\{ |\ip{a}{\lambda}| : \|\lambda^*\|\leq 1 \} \big\} \\
&= \max\big\{ \|a\|, \sup\{ |\ip{a^*}{\lambda}| : \|\lambda\|\leq 1 \} \big\}
= \max\big( \|a\|, \|a^*\| \big).
\end{align*}
Hence $\|a^*\|_0^* = \|a\|_0^*$, as required.
\end{proof}

Following \cite{Dales}, we shall say that a \emph{Banach $*$-algebra}
is a Banach algebra with an isometric involution.  We now know that
there is no loss of generality to talk about dual Banach $*$-algebras
as long as the involution is weak$^*$-continuous, which in light of
Proposition~\ref{star_has_rep} seems necessary for our purposes.
The next theorem shows that we can always embed a dual Banach algebra
with involution into a dual Banach algebra with weak$^*$-continuous
involution.  We remind the reader of Proposition~\ref{wap_dba}.

\begin{theorem}
Let $(\mc A,\mc A_*)$ be a dual Banach algebra with a continuous
involution.  Then $\wap(\mc A')'$ admits a weak$^*$-continuous
involution such that the canonical map $\mc A\rightarrow\wap(\mc A')'$
becomes a $*$-homomorphism.
\end{theorem}
\begin{proof}
By Proposition~\ref{wap_dba}, the first Arens product on $\mc A''$
drops to a well-defined product on $\wap(\mc A')'$ turning
$\wap(\mc A')'$ into a dual Banach algebra.  As shown in
\cite{DawLau}, it follows from Grothendieck's double limit criterion
that $\wap(\mc A')$ is a self-adjoint subspace of $\mc A'$;
in particular, this means that for $\mu\in\mc A_*$, certainly
$\mu^*\in\wap(\mc A')$.

We define an involution on $\mc A''$ by setting
\[ \ip{\Phi^*}{\lambda} = \overline{ \ip{\Phi}{\lambda^*} }
\qquad (\Phi\in\mc A'', \lambda\in\mc A'). \]
We then define an involution on $\wap(\mc A')'$ by setting
\[ (\Phi + \wap(\mc A')^\perp)^* = \Phi^* + \wap(\mc A')^\perp
\qquad (\Phi\in\mc A''), \]
which is well-defined, as $\wap(\mc A')$ is self-adjoint.
We may check that, for $a\in\mc A$, $\lambda\in\wap(\mc A')$ and
$\Phi\in\mc A''$, we have $(a\cdot\lambda)^* = \lambda^*\cdot a^*$
and $(\lambda\cdot\Phi)^* = \Phi^*\cdot\lambda$.
We then see that for $\Phi,\Psi\in\mc A''$ and $\lambda\in\wap(\mc A')$,
\begin{align*} \ip{(\Phi\aone\Psi)^*}{\lambda}
&= \overline{ \ip{\Phi}{\Psi\cdot\lambda^*} }
= \ip{\Phi^*}{(\Psi\cdot\lambda^*)^*}
= \ip{\Phi^*}{\lambda\cdot\Psi^*} \\
&= \ip{\Psi^*\atwo\Phi^*}{\lambda}
= \ip{\Psi^*\aone\Phi^*}{\lambda}, \end{align*}
by Lemma~\ref{wap_char}.  Thus $\wap(\mc A')'$ has a continuous
involution which, by definition, is weak$^*$-continuous, and
extends the involution on $\mc A$.
\end{proof}

Instead of using $\wap(\mc A')$ in the above construction, we
could instead have used $X = \mc A_* + \mc A_*^* \subseteq\wap(\mc A')$,
which has the advantage that if the involution on $\mc A$ is already
weak$^*$-continuous, then $X' = \mc A$.

\section{Connes-amenability and injectivity}

We shall show below that if $\mc A$ is Connes-amenable, then
$\mc A$ is unital.  In fact, a stronger result holds.

\begin{proposition}
Let $(\mc A,\mc A_*)$ be a dual Banach algebra such that
$(\mc A^\flat, \mc A_*^\flat)$ is Connes-amenable.
Then $\mc A$ is unital.
\end{proposition}
\begin{proof}
Let $E=\mc A_*$ as a Banach space, and for $\mu\in\mc A_*$,
write $\hat\mu$ for the canonical image of $\mu$ in $E$ (and
similarly for elements of $\mc A$ in $E'$).  Then turn $E$ into
an $\mc A^\flat$-bimodule by setting
\[ (a+\alpha)\cdot\hat\mu = \widehat{a\cdot\mu} + \alpha\hat\mu,\quad
\hat\mu\cdot (a+\alpha) = \alpha\hat\mu \qquad
(a+\alpha\in\mc A^\flat, \hat\mu\in E). \]
Then we claim that $E'$ becomes a normal $\mc A^\flat$-bimodule;
for example, if $a_i+\alpha_i\rightarrow a+\alpha$ weak$^*$ in
$\mc A^\flat$ then, for $\hat\mu\in E$ and $\hat b\in E'$,
\begin{align*}
\lim_i \ip{\hat b\cdot(a_i+\alpha_i)}{\hat\mu}
= \lim_i \ip{\hat b}{\widehat{a_i\cdot\mu} + \alpha_i\hat\mu}
&= \lim_i \ip{ba_i + \alpha_i b}{\mu} \\
= \ip{ba +\alpha b}{\mu}
&= \ip{\hat b \cdot (a+\alpha)}{\hat\mu},
\end{align*}
so that $\hat b\cdot(a_i+\alpha_i) \rightarrow \hat b\cdot
(a+\alpha)$ weak$^*$ in $E'$, as required.

Now let $d:\mc A^\flat\rightarrow E'$ be defined by $d(a+\alpha)
= \hat a$.  Then
\[ d( (a+\alpha)(b+\beta) ) = \widehat{ab} + \alpha\hat b +\beta\hat a
= \hat a \cdot (b+\beta) + (a+\alpha)\cdot\hat b, \]
so that $d$ is a derivation.  Clearly $d$ is weak$^*$-continuous,
so that as $(\mc A^\flat, \mc A_*^\flat)$ is Connes-amenable,
$d$ is inner, so that for some $e\in\mc A$, we have that
\[ \hat a = d(a+\alpha) = \hat e\cdot(a+\alpha) - 
(a+\alpha)\cdot \hat e = \widehat{ea} \qquad (a\in\mc A). \]
Thus $\mc A$ has a left identity, and in an analogous manner,
$\mc A$ has a right identity, so that $\mc A$ is unital.
\end{proof}

Thus the naive unitisation is useless as far as Connes-amenability
is concerned.

Let $E$ be a Banach space, and let $\mc A\subseteq\mc B(E)$ be
a subset.  We define the \emph{commutant} of $\mc A$ to be
\[ \mc A^c = \{ T\in\mc B(E) : TS=ST \ (S\in\mc A) \}, \]
so that $\mc A^c$ is a closed subalgebra of $\mc B(E)$.  We
then define $\mc A^{cc} = (\mc A^c)^c$, and see that $\mc A
\subseteq \mc A^{cc}$.

\begin{definition}
Let $E$ be a Banach space, and let $\mc A\subseteq\mc B(E)$ be
a subalgebra.  A \emph{quasi-expectation} for $\mc A$ is
a projection $\mc Q:\mc B(E)\rightarrow\mc A^c$ such that
$\mc Q(cTd) = c \mc Q(T) d$ for $c,d\in\mc A^c$ and $T\in\mc B(E)$.
\end{definition}

\begin{proposition}\label{ca_inj}
Let $\mc A$ be a Banach algebra, and let $\pi:\mc A\rightarrow\mc B(E)$
be a homomorphism with $E$ a reflexive Banach space.  Suppose that either:
\begin{enumerate}
\item $\mc A$ is amenable;
\item $\mc A$ is a dual Banach algebra, $\pi$ is weak$^*$-continuous,
   and $\mc A$ is Connes-amenable.
\end{enumerate}
Then  there exists a quasi-expectation
$\mc Q:\mc B(E) \rightarrow \pi(\mc A)^c$.
\end{proposition}
\begin{proof}
We may either translate, almost verbatim, the proof of
\cite[Theorem~3]{BP}, or else look at \cite[Theorem~4.4.11]{RundeBook}.
\end{proof}

Very similar ideas to the above are considered by Corach and Gal\'e
in \cite{CG}, and in particular Section~3 of that paper, where
they ask if the existance of a quasi-expectation is equivalent to
some form of amenability.  We shall now answer this question in
the affirmative, showing that quasi-expectations and Connes-amenability
are intimately linked (as is true in the von Neumann Algebra case: see \cite{BP}).

\begin{definition}
Let $(\mc A,\mc A_*)$ be a dual Banach algebra, and
let $E$ be a Banach $\mc A$-bimodule.  Then $x\in \sigma WC(E)$ if
and only if the maps $\mc A\rightarrow E$,
\[ a \mapsto \begin{cases} a\cdot x, \\ x\cdot a, \end{cases} \]
are $\sigma(\mc A,\mc A_*) - \sigma(E,E')$ continuous.
\noproof
\end{definition}

It is clear that $\sigma WC(E)$ is a closed submodule of $E$.
The $\mc A$-bimodule homomorphism $\Delta_{\mc A}$ has adjoint
$\Delta'_{\mc A}:\mc A' \rightarrow (\mc A\proten\mc A)'$.
In \cite[Corollary~4.6]{Runde2} it is shown that
$\Delta'_{\mc A}(\mc A_*) \subseteq \sigma WC( (\mc A\proten\mc A)' )$.
Consequently, we can view $\Delta_{\mc A}'$ as a map
$\mc A_* \rightarrow \sigma WC( (\mc A\proten\mc A)' )$, denoted by
$\tilde\Delta_{\mc A}$, and hence we have a map $\tilde\Delta_{\mc A}':
\sigma WC( (\mc A\proten\mc A)' )' \rightarrow \mc A_*'
= \mc A$.  The following is \cite[Theorem~4.8]{Runde2}.

\begin{theorem}\label{Runde_Thm}
Let $\mc A$ be a dual Banach algebra with predual $\mc A_*$.
Then the following are equivalent:
\begin{enumerate}
\item $\mc A$ is Connes-amenable;
\item $\mc A$ has a \emph{$\sigma WC$-virtual diagonal}, which
is $M\in \sigma WC( (\mc A\proten\mc A)' )'$ such that
$a\cdot M = M\cdot a$ and $a \tilde \Delta_{\mc A}'(M) = a$ for
each $a\in\mc A$.
\end{enumerate}
\end{theorem}

We can identify $\sigma WC( (\mc A\proten\mc A)' )'$ in a
slightly more concrete way.  First of all, recall that
$(\mc A\proten\mc A)' = \mc B(\mc A,\mc A')$, where we choose the
convention that
\[ \ip{T}{a\otimes b} = \ip{T(b)}{a} \qquad (a\otimes b\in\mc A
\proten\mc A, T\in\mc B(\mc A,\mc A')). \]
Then, for $\mu\in\mc A_*$, we identify $\tilde\Delta_{\mc A}(\mu)$
with the map $b\mapsto b\cdot\mu$.
The following is \cite[Proposition~3.2]{Connespp}.

\begin{proposition}\label{sigmawc_cond}
Let $(\mc A,\mc A_*)$ be a dual Banach algebra.
For $T\in\mc B(\mc A,\mc A') = (\mc A\proten\mc A)'$, define
maps $\phi_r, \phi_l : \mc A\proten\mc A\rightarrow\mc A'$ by
\[ \phi_r(a\otimes b) = T'\kappa_{\mc A}(a) \cdot b,
\quad \phi_l(a\otimes b) = a \cdot T(b)
\qquad (a\otimes b\in\mc A\proten\mc A). \]
Then $T\in\sigma WC(\mc B(\mc A,\mc A'))$ if and only if
$\phi_r$ and $\phi_l$ are weakly-compact and have ranges
contained in $\tobidual{\mc A_*}$.
\end{proposition}

As the unit ball of $\mc A\proten\mc A$ is the closure of the convex
hull of $\{ a\otimes b : \|a\| = \|b\|\leq 1\}$, we see that, for example,
$\phi_l$ is weakly-compact if and only if the set $\{a\cdot T(b) :
\|a\| = \|b\|\leq 1\}$ is relatively weakly (sequentially) compact.

We shall now prove a representation result for maps in
$\sigma WC( \mc B(\mc A,\mc A') )$.  Firstly, we again
apply some interpolation space theory.

\begin{proposition}\label{sigmawc_int}
Let $(\mc A,\mc A_*)$ be a unital dual Banach algebra, and let
$T\in\sigma WC(\mc B(\mc A,\mc A'))$.  Then there exists an absolute
constant $K>0$, a Banach left $\mc A$-module $E$ and a Banach right
$\mc A$-module $F$ such that $E$ and $F$ are normal and reflexive,
and for some unit vectors $\mu_0\in E'$ and $\lambda_0\in F'$,
we have that
\[ |\ip{T(b)}{a}| \leq K \|T\| \| \mu_0\cdot a \|_{E'}
\|b\cdot\lambda_0\|_{F'} \qquad (a,b\in\mc A). \]
Furthermore, $\mu_0\cdot\mc A$ is dense in $E'$ and
$\mc A\cdot\lambda_0$ is dense in $F'$.
\end{proposition}
\begin{proof}
We may suppose that $\|T\|=1$.  Form $\phi_l$ using $T$, and for
$n\geq1$ define a norm on $\mc A_*$ by
\[ \|\mu\|_n = \inf\big\{ 2^{-n/2}\|\tau\| +
2^{n/2}\|\mu-\phi_l(\tau)\| : \tau\in\mc A\proten\mc A \big\}
\qquad (\mu\in\mc A_*). \]
As in the proof of Theorem~\ref{int_thm}, we may check that $\|\cdot\|_n$
is an equivalent norm on $\mc A_*$, and that we may define
\[ E = \Big\{ \mu\in\mc A_* : \|\mu\|_E := \Big(
\sum_{n\geq1} \|\mu\|_n^2 \Big)^{1/2} < \infty \Big\}. \]
Then $E$ is a Banach space, $\phi_l(\mc A\proten\mc A)$ is dense in
$E$, and for $\tau\in\mc A\proten\mc A$, we have that
$\| \phi_l(\tau) \| \leq \|\phi_l(\tau)\|_E \leq \| \tau \|_\pi$.
Furthermore, $E$ is a Banach left $\mc A$-module as $\phi_l$ is a
left $\mc A$-module homomorphism.
Again, it follows from general interpolation space results that
$E$ is reflexive, as $\phi_l$ is weakly-compact.

Again, let $\iota_E:E\rightarrow\mc A_*$ be the inclusion map,
so that $\iota_E':\mc A\rightarrow E'$ has dense range.  Let $\mu\in E$,
$\iota_E'(a)\in E'$, and let $(b_\alpha)$ be a net in $\mc A$
which converges weak$^*$ to $b$.  Then
\[ \ip{\iota_E'(a)}{b_\alpha\cdot\mu} =
\ip{ab_\alpha}{\mu} \rightarrow \ip{ab}{\mu} = \ip{\iota_E'(a)}{b\cdot\mu}, \]
so we see that the map $\mc A\rightarrow E, b\mapsto b\cdot \mu$ is
weak$^*$-continuous, that is, $E$ is normal.  For each $n\geq1$,
let $\|\cdot\|_n^*$ be the dual norm to $\|\cdot\|_n$, defined on
$\mc A$.  For $a\in\mc A$, we have that
\begin{align*} \|a\|_n^* &= \sup\big\{ |\ip{a}{\mu + \phi_l(\tau)}| :
   2^{-n/2}\|\tau\| + 2^{n/2}\|\mu\| \leq 1 \big\} \\
&= \sup\big\{ 2^{-n/2}|\ip{a}{\mu}| + 2^{n/2}|\ip{\phi_l'(a)}{\tau}| :
   \|\tau\| + \|\mu\| \leq 1 \big\} \\
&= \max\big( 2^{-n/2}\|a\| , 2^{n/2}\|\phi_l'(a)\| \big).
\end{align*}
As $E$ isometrically embeds into the $\ell^2$-direct sum of the spaces
$(\mc A_*,\|\cdot\|_n)_{n\geq 1}$, we see that $E'$ is isometrically
a quotient of the $\ell^2$-direct sum of the spaces
$(\mc A,\|\cdot\|_n^*)_{n\geq 1}$.  We hence see that, if we drop
the map $\iota_E'$ and identify $\mc A$ with a dense subspace of
$E'$, then for $a\in\mc A$,
\[ \|a\|_{E'} = \inf\Big\{ \Big( \sum_{n\geq 1}
\max\big( 2^{-n/2}\|a_n\| , 2^{n/2}\|\phi_l'(a_n)\| \big)^2 \Big)^{1/2}
: a = \sum_{n\geq1} a_n \Big\}. \]
From \cite[Section~1.5, Proposition~1]{BB} it follows that there
is an absolute constant $K>0$ such that if we define
\[ \|a\|_1 = \inf\Big\{ \Big( \sum_{n\geq 1}
   \big( 2^{-n/2}\|a-a_n\| + 2^{n/2}\|\phi_l'(a_n)\| \big)^2 \Big)^{1/2}
   : (a_n)\subseteq\mc A \Big\}, \]
then $K^{-1} \|\cdot\|_1 \leq \|\cdot\|_{E'} \leq K \|\cdot\|_1$.

We analogously use $\phi_r$ to form a reflexive, normal, Banach 
right $\mc A$-module $F$, and we find a norm $\|\cdot\|_2$ on $\mc A$
such that $K^{-1} \|\cdot\|_2 \leq \|\cdot\|_{F'} \leq K \|\cdot\|_2$.
Notice that for $a,b\in\mc A$, $\ip{T(b)}{a} =
\ip{\phi_l'(a)(b)}{e_{\mc A}}$ and
\begin{align*} \|\phi_l'(a)(b)\| &= \sup\{ |\ip{T(b)}{ac}| : \|c\|\leq1 \} \\
&= \sup\{ |\ip{\phi_r'(b)(e_{\mc A})}{ac}| : \|c\|\leq1 \}
\leq \| \phi_r'(b) \| \|a\|. \end{align*}
For $a,b\in\mc A$ and $\epsilon>0$, choose $(a_n),(b_n)\subseteq\mc A$
such that $a=\sum_{n\geq 1} a_n$ and
\begin{align*}
\Big( \sum_{n\geq 1} \max\big( 2^{-n/2}\|a_n\| ,
2^{n/2}\|\phi_l'(a_n)\| \big)^2 \Big)^{1/2} \leq \|a\|_{E'} + \epsilon, \\
\Big( \sum_{n\geq 1} \big( 2^{-n/2}\|b-b_n\| +
2^{n/2}\|\phi_r'(b_n)\| \big)^2 \Big)^{1/2} \leq \|b\|_2 + \epsilon.
\end{align*}
We then use the Cauchy-Schwarz inequality to see that
\begin{align*}
|\ip{T(b)}{a}| &\leq \|\phi_l'(a)(b)\| \leq \sum_{n\geq1} \|\phi_l'(a_n)(b)\| \\
&\leq \sum_{n\geq1} \|\phi_l'(a_n)(b_n-b)\| + \|\phi_l'(a_n)(b_n)\| \\
&\leq \sum_{n\geq1} 2^{n/2}\|\phi_l'(a_n)\| 2^{-n/2}\|b_n-b\| +
   \sum_{n\geq 1} 2^{n/2}\|\phi_r'(b_n)\| 2^{-n/2}\|a_n\| \\
&\leq 2(\|a\|_{E'}+\epsilon)(\|b\|_2+\epsilon).
\end{align*}
As $\epsilon>0$ was arbitrary, we see that $|\ip{T(b)}{a}| \leq
2K \|a\|_{E'} \|b\|_{F'}$.  We may hence set $\mu_0 = \iota_E'(e_{\mc A})$
and $\lambda_0 = \iota_F'(e_{\mc A})$, so that as $\iota_E'$ has
dense range, $\mu_0\cdot\mc A$ is dense in $E'$, and similarly
for $\mc A\cdot\lambda_0$.  It remains to check that
these are unit vectors.  However, as $\iota'$ is norm-decreasing,
$\| \mu_0 \|_{E'} \leq 1$, while conversely
\begin{align*} \|\iota_E'(e_{\mc A})\| &=
\sup\{ |\ip{e_{\mc A}}{\phi_l(\tau)}| : \| \phi_l(\tau) \|_E\leq 1\} \\
&\geq \sup\{ |\ip{e_{\mc A}}{\phi_l(a\otimes b)}| :
   \| a\otimes b\|_\pi = \|a\| \|b\| \leq 1 \} \\
&= \sup\{ |\ip{a}{T(b)}| : \|a| = \|b\| = 1 \} = \|T\| = 1.
\end{align*}
Similarly, $\|\lambda_0\|_{F'}=1$, and the proof is complete.
\end{proof}

Notice that if $T = \tilde\Delta_{\mc A}(\mu)$ for some $\mu\in
\mc A_*$, then $\phi_l(a\otimes b) = a\cdot (b\cdot\mu) = ab\cdot\mu$.
It hence follows that the space $E$ constructed in the above
proof is equal to the space constructed in Theorem~\ref{int_thm}.

\begin{theorem}\label{sigmawc_rep}
Let $(\mc A,\mc A_*)$ be a unital dual Banach algebra, and let
$T\in\mc B(\mc A,\mc A')$.  Then the following are equivalent:
\begin{enumerate}
\item $T \in \sigma WC( \mc B(\mc A,\mc A') )$;
\item there exists a normal reflexive Banach left $\mc A$-module
   $E$, $x\in E$, $\mu\in E'$ and $S\in\mc B(E)$ such that
   $\ip{T(b)}{a} = \ip{\mu\cdot a}{S(b\cdot x)}$ for each $a,b\in\mc A$.
\end{enumerate}
Furthermore, in this case, there is an absolute constant $K>0$
such that we may choose $E$, $x$, $\mu$ and $S$ with
$\|x\| \|\mu\| \|S\| \leq K \|T\|$.
\end{theorem}
\begin{proof}
If (2) holds, then let $\pi:\mc A\rightarrow\mc B(E)$ be the
associated representation, and let $\pi_*: E'\proten E\rightarrow
\mc A'$ be the restriction of $\pi'$ to $E'\proten E$.  As $E$
is normal, the representation is weak$^*$-continuous, and so $\pi_*$
maps into $\mc A_*$, and $\pi_*' = \pi$.
For $a,b\in\mc A$, we have that
\begin{align*} \ip{T(b)}{a} &= \ip{\mu}{\pi(a)S\pi(b)(x)}
= \ip{a}{\pi_*\big( S\pi(b)(x) \otimes \mu \big)} \\
&= \ip{b}{\pi_*\big( x \otimes S'\pi(a)'(\mu) \big)}, \end{align*}
so that $T(\mc A)\subseteq \mc A_*$ and $T'\kappa_{\mc A}(\mc A)
\subseteq\mc A_*$.

Let $(a_n)$ and $(b_n)$ be bounded sequences in $\mc A$, and let $x\in E$.
As $E$ is reflexive, the unit ball of $E$ is weakly sequentially
compact, and so, by moving to subsequences if necessary, we may
suppose that for some $y\in E$,
\[ \ip{\mu}{\pi(a_n)S\pi(b_n)(x)} \rightarrow \ip{\mu}{y}
\qquad (\mu\in E'). \]
Then, for $c\in\mc A$ and $\mu\in E'$,
\begin{align*} \lim_{n\rightarrow\infty}&\ip{a_n\cdot T(b_n)}{c} =
\lim_{n\rightarrow\infty} \ip{\mu}{\pi(ca_n)S\pi(b_n)(x)} \\
&= \lim_{n\rightarrow\infty} \ip{\pi(c)'(\mu)}{\pi(a_n)S\pi(b_n)(x)}
= \ip{\pi(c)'(\mu)}{y} = \ip{c}{\pi_*(y\otimes\mu)}.
\end{align*}
Hence, combined with the comments after Proposition~\ref{sigmawc_cond},
we see that $\phi_l$ (as defined using $T$) is weakly-compact.
A similar calculation shows that $\phi_r$ is weakly-compact,
and so we conclude that $T\in\sigma WC( \mc B(\mc A,\mc A'))$
as required.

Conversely, let $T\in\sigma WC(\mc B(\mc A,\mc A'))$, and form
the spaces $E$ and $F$ using Proposition~\ref{sigmawc_int}.
We have that $\mu_0\cdot\mc A$ is dense in $E'$ and that
$\mc A\cdot\lambda_0$ is dense in $F'$.
Then define $R\in (E'\proten F')' = \mc B(F',E)$ by
\[ \ip{R}{\mu_0\cdot a \otimes b\cdot\lambda_0} = \ip{T(b)}{a}
\qquad (a,b\in\mc A), \]
so that $\|R\| \leq K\|T\|$.
Let $G = E \oplus_2 F'$ (that is, the norm on $G$ is
$\| (e,f) \| = (\|e\|^2 + \|f\|^2)^{1/2}$ for $e\in E$ and $f\in F$)
so that $G' = E' \oplus_2 F$, and $G$ is a normal, reflexive,
Banach left $\mc A$-module.  Define $S\in\mc B(G)$ by
\[ S(z,b\cdot\lambda_0) = (R(b\cdot\lambda_0), 0)
\qquad (z\in E, b\cdot\lambda_0\in F'), \]
and let $x = (0,\lambda_0)\in G$, $\mu = (\mu_0,0)\in G'$.
Then, for $a,b\in\mc A$,
\[ \ip{\mu\cdot a}{S(b\cdot x)} =
\ip{(\mu_0\cdot a,0)}{R(b\cdot\lambda_0,0)}
= \ip{R}{\mu_0\cdot a \otimes b\cdot\lambda_0}
= \ip{T(b)}{a}, \]
as required.
\end{proof}

Let $\mc A$ be a Banach algebra and $E$ be a left $\mc A$-module.
Then we write ${_{\mc A}}\mc B(E)$ for the collection of \emph{left
$\mc A$-module homomorphisms}, that is, maps $T\in\mc B(E)$ such that
$T(a\cdot x)=a\cdot T(x)$ for $a\in\mc A$ and $x\in E$.  Similarly,
we define $\mc B_{\mc A}(E)$ and $\bimodhomo{\mc A}{E}$ to be the
collection of \emph{right $\mc A$-module homomorphisms} and
\emph{$\mc A$-bimodule homomorphisms}, respectively.

Suppose now that $\mc A$ is a closed subalgebra of $\mc B(E)$ for
some reflexive Banach space $E$.  Then $\mc B(E)$ becomes a
Banach $\mc A$-bimodule and a Banach $\mc A^c$-bimodule in the
obvious way.  We turn $\bimodhomo{\mc A^c}{\mc B(E)}$ into
a Banach $\mc A$-module by setting
\[ (a\cdot\mc S)(T) = a \mc S(T), \quad
(\mc S\cdot a)(T) = \mc S(T)a \qquad (a\in\mc A, T\in\mc B(E) ), \]
for $\mc S \in \bimodhomo{\mc A^c}{\mc B(E)}$.
Notice that $\mc B(\mc B(E))$ is a dual Banach space with
predual $\mc B(E) \proten (E \proten E')$.  Let $X
\subseteq \mc B(E) \proten (E \proten E')$ be the closure
of the linear span of
\[ \left\{ \begin{matrix} cT \otimes x\otimes \mu - T\otimes x\otimes c'(\mu), \\
Tc\otimes x\otimes\mu - T\otimes c(x)\otimes\mu \end{matrix}
: c\in\mc A^c, T\in\mc B(E), x\in E,\mu\in E' \right\}. \]
Then, for example, if $\mc S\in\mc B(\mc B(E))$ is such that
$\ip{\mc S}{cT\otimes x\otimes\mu - T\otimes x\otimes c'(\mu)}=0$
for each $c\in\mc A^c, T\in\mc B(E), x\in E$ and $\mu\in E'$, then
$\mc S(cT) = c\mc S(T)$.  We hence see that $X^\perp =
\bimodhomo{\mc A^c}{\mc B(E)}$, so that $\bimodhomo{\mc A^c}{\mc B(E)}$
has the predual $\mc B(E)\proten E\proten E' / X$.

Define $\theta:\mc A\proten\mc A\rightarrow\bimodhomo{\mc A^c}{\mc B(E)}$
by
\[ \theta(a\otimes b)(T) = aTb \qquad (a,b\in\mc A, T\in\mc B(E)), \]
so that $\theta$ is an $\mc A$-bimodule homomorphism.  We then
define $\psi: \mc B(E)\proten E\proten E' / X \rightarrow
\mc B(\mc A,\mc A')$ by, for $a,b\in\mc A, x\in E, \mu\in E'$ and
$T\in\mc B(E)$,
\[ \ip{\psi(T\otimes x\otimes\mu + X)}{a\otimes b} = 
\ip{\mu}{\theta(a\otimes b)(T)(x)} = \ip{\mu}{aTb(x)}. \]
A simple check shows that this is well-defined, and that $\|\psi\|\leq1$.
We turn $\mc B(E) \proten E\proten E'$ into a Banach $\mc A$-bimodule
by setting
\[ \begin{matrix} a\cdot (T\otimes x\otimes\mu) =
T\otimes a(x)\otimes\mu, \\ (T\otimes x\otimes\mu)\cdot a =
T\otimes x\otimes a'(\mu) \end{matrix}
\qquad (a\in\mc A, T\in\mc B(E), x\in E, \mu\in E'). \]
Then $X$ is a sub-$\mc A$-bimodule, and this module action
agrees with the module action already defined on
$\bimodhomo{\mc A^c}{\mc B(E)}$.  We may verify that $\psi$ is
an $\mc A$-bimodule homomorphism.

We now aim to construct a Banach space $E$ such that $\psi$ is
a bijection (onto a suitable closed subspace of $\mc B(\mc A,\mc A')$)
for this $E$.

\begin{definition}
For a Banach left $\mc A$-module $E$, and $x\in E$, let
$\overline{\mc A\cdot x}$ be the closure of $\{ a\cdot x: a\in\mc A\}$,
so that $\overline{\mc A\cdot x}$ is a closed submodule of $E$.
Similarly, for $\mu\in E'$, define $\overline{\mu\cdot\mc A}$.
We then say that $E$ is \emph{cyclic} if, for some $x_0\in E$,
we have that $\overline{\mc A\cdot x_0} = E$, and a similar definition
holds for $E'$.
\end{definition}

For Banach spaces $E$ and $F$, we let $l^2(E\oplus F) = E\oplus_2 F$
be the direct sum of $E$ and $F$ with the norm $\| (e,f) \| =
( \|e\|^2 + \|f\|^2 )^{1/2}$ for $e\in E, f\in F$.  When $E$ and
$F$ are reflexive, normal, Banach left $\mc A$-modules, it is clear
that $E\oplus_2 F$ is also.  We similarly define
$l^2\big(\bigoplus_\alpha E_\alpha\big)$, where $(E_\alpha)$ is a
family of Banach spaces.

The following lemma is a technical result.  It would be easier to
define $E$ to be the $l^2$-direct-sum of \emph{all} reflexive, normal,
Banach left $\mc A$-modules, but this collection is not in general
a set.

\begin{lemma}\label{set_theory_ok}
Let $\mc E$ be a set of reflexive, normal, cyclic, Banach left
$\mc A$-modules.  There exists a reflexive, normal, Banach left
$\mc A$-module $E$ such that:
\begin{enumerate}
\item each member of $\mc E$ is isometrically isomorphic to a
   $1$-complemented submodule of $E$;
\item for $x_1,x_2\in E$ and $\mu_1,\mu_2\in E'$, if we let
   $X(x_1,x_2) = \overline{ \{ a\cdot(x_1,x_2) : a\in\mc A \}}
   \subseteq \overline{\mc A\cdot x_1} \oplus_2
   \overline{\mc A\cdot x_2}$, and we let $Y(\mu_1,\mu_2) =
   \overline{ \{ (\mu_1,\mu_2)\cdot a : a\in\mc A \}} \subseteq
   \overline{\mu_1\cdot\mc A} \oplus_2 \overline{\mu_2\cdot\mc A}$,
   then $X(x_1,x_2) \oplus_2 Y(\mu_1,\mu_2)'$ is isometrically
   isomorphic to a $1$-complemented submodule of $E$;
\end{enumerate}
\end{lemma}
\begin{proof}
Let $\mc E_0 = \mc E$, and then use transfinite induction
to define, for an ordinal $\alpha$, $\mc E_\alpha$ as follows.
If $\alpha$ is a limit ordinal, we let $\mc E_\alpha = \bigcup_{\lambda
<\alpha} \mc E_\lambda$.  Otherwise, let $E_\alpha = 
l^2(\bigoplus_{E\in\mc E_\alpha} E)$ and for $x_1,x_2\in E_\alpha$ and
$\mu_1,\mu_2\in E_\alpha'$, form $X(x_1,x_2)$ and $Y(\mu_1,\mu_2)$
as above.  Then let $\mc E_{\alpha+1}$ be $\mc E_{\alpha}$ unioned
with the collection of all such spaces $X(x_1,x_2) \oplus_2
Y(\mu_1,\mu_2)'$.  Notice that each member of $\mc E_\alpha$ is
canonically a reflexive, normal, Banach $\mc A$-module.
Notice that $X(x_1,x_2)$ is always
a cyclic module, while $Y(\mu_1,\mu_2)'$ is the dual of a cyclic
(right) module.

We then let $E = l^2(\bigoplus_{E\in\mc E_{\aleph_1}} E)$ and give
$E$ the obvious left $\mc A$-module structure.  Then each member of
$\mc E$ is a $1$-complemented subspace of $E$.  Indeed, we may view
$E_\alpha = l^2(\bigoplus_{E\in\mc E_\alpha} E)$ and $E_\alpha'$
as submodules of $E$ and $E'$, respectively, for each $\alpha<\aleph_1$.
For notational convenience, let $\mc E_{\aleph_1} = \{ E_i : i\in I \}$
for some indexing set $I$.  Then let $x_1,x_2\in E$ and
$\mu_1,\mu_2\in E'$, so that, for $k=1,2$, $x_k = (x^{(k)}_i)_{i\in I}$
with $\|x_k\| = \big( \sum_{i\in I} \|x^{(k)}_i\|^2 \big)^{1/2}$,
and similarly for $\mu_k$.  As $\aleph_1$ is the first uncountable
ordinal, we see that for some $\alpha<\aleph_1$, for each $k=1,2$
we have that $x^{(k)}_i\not=0$ or $\mu^{(k)}_i\not=0$ implies that
$E_i\in\mc E_\alpha$.  Hence $x_1,x_2\in E_\alpha$ and $\mu_1,\mu_2
\in E_\alpha'$, so that by construction, $X(x_1,x_2) \oplus_2
Y(\mu_1,\mu_2)'$ is a $1$-complemented submodule of
$E_{\alpha+1} \subseteq E$ as required.
\end{proof}

\begin{theorem}\label{rep_psi_bi}
Let $(\mc A,\mc A_*)$ be a unital dual Banach algebra.
There exists an isometric, weak$^*$-continuous representation $\pi:\mc A
\rightarrow \mc B(E)$ such that $\psi$ (as associated with $\pi$)
maps into $\sigma WC( \mc B(\mc A,\mc A'))$ and is a bijection.
\end{theorem}
\begin{proof}
By Theorem~\ref{sigmawc_rep}, we see that $\psi$ maps into
$\sigma WC( \mc B(\mc A,\mc A'))$ for any isometric weak$^*$-continuous
representation $\pi:\mc A\rightarrow\mc B(E)$.  Let $\mc E$ be
the collection of Banach spaces constructed in Theorem~\ref{sigmawc_rep}
for each norm-one member of $\sigma WC(\mc B(\mc A,\mc A'))$.  Then
let $E$ be the Banach space given by Lemma~\ref{set_theory_ok},
so that it is clear that $\psi$ is surjective, by Theorem~\ref{sigmawc_rep}.
As such, $\psi'$ is an isomorphism onto its range.  We shall now show
that $\psi'$ is surjective, which will complete the proof.

Fix $\mc S\in \bimodhomo{\mc A^c}{\mc B(E)}$, and define $M\in
\sigma WC(\mc B(\mc A,\mc A'))'$ in the following way.  For each
$T\in\sigma WC(\mc B(\mc A,\mc A'))$, let $x\in E$, $\mu\in E'$ and
$S\in\mc B(E)$ be such that $\ip{T}{a\otimes b} =
\ip{\mu\cdot a}{S(b\cdot x)}$ for each $a,b\in\mc A$.  Then define
$\ip{M}{T} = \ip{\mu}{\mc S(S)(x)}$.  Suppose that this \emph{is}
well-defined.  Then, for each $x\in E, \mu\in E'$ and $S\in\mc B(E)$,
\begin{align*} \ip{\psi'(M)}{S\otimes x\otimes\mu+X}
&= \ip{M}{\psi(S\otimes x\otimes \mu+X)} = \ip{\mu}{\mc S(S)(x)} \\
&= \ip{\mc S}{S\otimes x\otimes\mu+X}, \end{align*}
so that $\psi'(M) = \mc S$ as required.

We shall now show that $M$ is well-defined, at least for our specific $E$.
Let $T\in\sigma WC(\mc B(\mc A,\mc A'))$, and suppose that, for $i=1,2$,
we have $x_i\in E$, $\mu_i\in E'$ and $S_i\in\mc B(E)$ such that
$\ip{T}{a\otimes b} = \ip{\mu_i\cdot a}{S_i(b\cdot x_i)}$ for $a,b\in\mc A$.
Pick $t\in (0,1)$ such that $t\|S_1\| = (1-t)\|S_2\|$.  A quick
calculation shows that then $K := t\|S_1\| = \|S_1\|\|S_2\|
(\|S_1\|+\|S_2\|)^{-1}$.
For each $a,b\in\mc A$, by the Cauchy-Schwarz inequality,
\begin{align*}
|\ip{T}{a\otimes b}| &= t|\ip{\mu_1\cdot a}{S_1(b\cdot x_1)}|
   + (1-t)|\ip{\mu_2\cdot a}{S_2(b\cdot x_2)}| \\
&\leq t\|S_1\| \|\mu_1\cdot a\| \|b\cdot x_1\|
   + (1-t)\|S_2\| \|\mu_2\cdot a\| \|b\cdot x_2\| \\
&\leq K \big( \|\mu_1\cdot a\|^2 + \|\mu_2\cdot a\|^2 \big)^{1/2}
   \big( \|b\cdot x_1\|^2 + \|b\cdot x_2\|^2 \big)^{1/2}.
\end{align*}

Let $F$ be the closure of $\mc A\cdot(x_1,x_2)$ in
$\overline{\mc A\cdot x_1} \oplus_2 \overline{\mc A\cdot x_2}$,
and let $G$ be the closure of $(\mu_1,\mu_2)\cdot\mc A$ in
$\overline{\mu_1\cdot\mc A} \oplus_2 \overline{\mu_2\cdot\mc A}$.
The above calculation allows us to define $R\in (G\proten F)'=
\mc B(F,G')$ by
\[ \ip{R}{(\mu_1\cdot a,\mu_2\cdot a) \otimes (b\cdot x_1,b\cdot x_2)}
= \ip{T}{a\otimes b} \qquad (a,b\in\mc A), \]
and we see that $\|R\|\leq K$.  Then set $H = G'\oplus_2 F$,
and let $P_F$ and $P_{G'}$ be the projections onto $F$ and
$G'$ respectively.  As $H'=G\oplus_2 F'$, let $P_G$ and
$P_{F'}$ be defined similarly, so that $P_G = P_{G'}'$ and
$P_{F'} = P_F'$.  Let $x_0=(0,(x_1,x_2))\in H$,
$\mu_0=((\mu_1,\mu_2),0)\in H'$,  and define
$S_0\in\mc B(H)$ by $S_0(g^*,f) = (R(f),0)$ for
$g^*\in G'$ and $f\in F$.  Then, for $a,b\in\mc A$,
\[ \ip{\mu_0\cdot a}{S_0(b\cdot x_0)} 
= \ip{(\mu_1\cdot a,\mu_2\cdot a)}{R(b\cdot x_1,b\cdot x_2)}
= \ip{T}{a\otimes b}. \]

By Lemma~\ref{set_theory_ok}, $H$ is a $1$-complemented
submodule of $E$, so we can find norm-decreasing left $\mc A$-module
homomorphisms $P:E\rightarrow H$ and $\iota:H\rightarrow E$
such that $P\iota = I_H$.
For $i=1,2$, we may define maps
$U_i \in {_{\mc A}}\mc B(H,\overline{\mc A\cdot x_i})$ and
$V_i \in \mc B_{\mc A}(H',\overline{\mu_i\cdot\mc A})$ by,
for $a\in\mc A, f^*\in F'$ and $g^*\in G'$,
\[ U_i(g^*,a\cdot(x_1,x_2)) = a\cdot x_i, \quad
V_i((\mu_1,\mu_2)\cdot a,f^*) = \mu_i\cdot a. \]
Then, by construction, $U_i$ and $V_i$ are norm-decreasing.
For $a,b\in\mc A$, we see that
\begin{align*}
\ip{T}{a\otimes b} &= \ip{\mu_i\cdot a}{S_i(b\cdot x_i)}
= \ip{\mu_0\cdot a}{V_i'S_iU_i(b\cdot x_0)} \\
&= \ip{\mu_0\cdot a}{P_{G'} V_i'S_iU_iP_F(b\cdot x_0)}
= \ip{\mu_0\cdot a}{P_{G'} S_0P_F(b\cdot x_0)}
\end{align*}
As $\mc A\cdot x_0$ is dense in $F$ and $\mu_0\cdot\mc A$ is
dense in $G$, we conclude that $P_{G'} V_i'S_iU_iP_F =
P_{G'} S_0P_F$.
As $\mc S \in \bimodhomo{\mc A^c}{\mc B(E)}$, we see that
\begin{align*}
\ip{\mu_i}{\mc S(S_i)(x_i)} &= 
\ip{V_iP_G(\mu_0)}{\mc S(S_i) U_iP_F(x_0)} \\
&= \ip{\mu_0}{P\iota P_{G'} V_i' \mc S(S_i) U_iP_F P\iota(x_0)} \\
&= \ip{\mu_0}{P \mc S(\iota P_{G'} V_i' S_i U_iP_F P)\iota(x_0)} \\
&= \ip{\mu_0}{P \mc S(\iota P_{G'} S_0 P_F P)\iota(x_0)}.
\end{align*}
We hence conclude that $\ip{\mu_i}{\mc S(S_i)(x_i)}$ has
the same value for $i=1$ as for $i=2$, and hence that $M$
is well-defined, as required.
\end{proof}

\begin{definition}
Let $(\mc A,\mc A_*)$ be a dual Banach algebra.  We say that
$\mc A$ is \emph{injective} if whenever $\pi:\mc A\rightarrow\mc B(E)$
is a weak$^*$-continuous, unital representation, there is a
quasi-expectation $\mc Q:\mc B(E)\rightarrow \pi(\mc A)^c$.
\end{definition}

\begin{theorem}\label{conn_inje}
Let $(\mc A,\mc A_*)$ be a unital dual Banach algebra.  Then $\mc A$ is
Connes-amenable if and only if $\mc A$ is injective.
\end{theorem}
\begin{proof}
We have already seen that when $\mc A$ is Connes-amenable, $\mc A$
is injective.  Conversely, consider the weak$^*$-continuous
representation $\pi:\mc A\rightarrow\mc B(E)$ constructed in 
Theorem~\ref{rep_psi_bi}, so that $\psi'$ is an isomorphism.
As $\mc A$ is injective, there exists a quasi-expectation
$\mc Q:\mc B(E)\rightarrow\pi(\mc A)^c$.  Let $M = (\psi')^{-1}(\mc Q)
\in \sigma WC((\mc A\proten\mc A)')'$.  As $\mc Q$ maps into
$\pi(\mc A)^c$, it follows that $a\cdot\mc Q = \mc Q\cdot a$ for
$a\in\mc A$, so that $a\cdot M = M\cdot a$.

As the unit ball of $\mc A\proten\mc A$ is weak$^*$-dense in
the unit ball of $\sigma WC((\mc A\proten\mc A)')'$, there exists
a bounded net $(\tau_\alpha)$ in $\mc A\proten\mc A$ such that
$M$ is the weak$^*$-limit of $(\tau_\alpha)$.  For each $\alpha$,
let $\tau_\alpha = \sum_{n\geq 1} a^{(\alpha)}_n \otimes
b^{(\alpha)}_n$.  For $x\in E$ and $\mu\in E'$, there exists
$\lambda\in\mc A_*$ such that $\ip{a}{\lambda} = \ip{\mu}{\pi(a)(x)}$
for $a\in\mc A$.  We then see that for $a,b\in\mc A$,
\[ \ip{ab}{\lambda} = \ip{\tilde\Delta_{\mc A}(\lambda)}{a\otimes b}
= \ip{\mu}{\pi(ab)(x)} = \ip{\mu}{\pi(a)I_E\pi(b)(\mu)}. \]
Then, from the proof of Theorem~\ref{rep_psi_bi}, we see that
\begin{align*}
\ip{\mu}{\mc Q(I_E)(x)} &= \ip{M}{\tilde\Delta_{\mc A}(\lambda)}
= \ip{\tilde\Delta_{\mc A}'(M)}{\lambda}
= \lim_\alpha \sum_{n\geq 1}
   \ip{\mu}{\pi(a^{(\alpha)}_nb^{(\alpha)}_n) (x)}.
\end{align*}
As $\mc Q$ is a projection onto $\pi(A)^c$ and $I_E\in\pi(A)^c$,
we see that $\mc Q(I_E) = I_E$, and so, as $x$ and $\mu$ were
arbitrary, we must have that $\lim_\alpha \sum_{n\geq1}
a^{(\alpha)}_nb^{(\alpha)}_n = e_{\mc A}$ in the weak$^*$-topology
on $\mc A$.  That is, $\tilde\Delta'_{\mc A}(M) = e_{\mc A}$,
showing that $M$ is a $\sigma WC$-virtual diagonal, which implies
that $\mc A$ is Connes-amenable, as required.
\end{proof}

There exists a rather strong decomposition theory for
weak$^*$-continuous homomorphisms between von Neumann algebras
(see \cite[Theorem~5.5]{tak}).  From this, it follows that if
$\mc A\subseteq\mc B(H)$ is a von Neumann algebra admitting
a quasi-expectation, and $\mc A$ is isomorphic to $\mc B\subseteq
\mc B(K)$, then $\mc B$ admits a quasi-expectation (see
\cite[Lemma~6.1.2]{RundeBook}).  Hence we need only look at one
representation for $\mc A$ to decide if $\mc A$ is Connes-amenable.

In contrast, it follows from \cite[Corollary~4.5]{Runde1} that
$\mc A = \mc B(\ell^p \oplus \ell^q)$ is not Connes-amenable when
$p,q\in (1,\infty)\setminus\{2\}$ are distinct, while trivially,
there is a quasi-expectation for $\mc A$ under the trivial
representation to $\mc B(\ell^p \oplus \ell^q)$.
Our theorem shows that there exists some reflexive Banach space
$E$ and some weak$^*$-continuous representation $\pi:
\mc A\rightarrow\mc B(E)$ such that $\mc A$ has no
quasi-expectation for $\mc B(E)$.  It would be interesting to
determine the Banach space properties of $E$.

There exists a more category-theoretic
definition of injectivity for von Neumann algebras (see
\cite[Chapter~XV, Section~1]{tak3}): that is, they are \emph{injective}
in the usual mapping sense, with respect to \emph{completely
positive maps}.  Does a similar definition hold for dual Banach algebras?

\begin{remark}
Let $(\mc A,\mc A_*)$ be a dual Banach algebra, and let
$\pi:\mc A\rightarrow\mc B(E)$ be the representation given by
Theorem~\ref{rep_psi_bi}.  Then suppose that $\mc B = \mc A^{cc}$
is Connes-amenable, so that there exists a quasi-projection
$\mc Q:\mc B\rightarrow(\mc A^{cc})^c = \mc A^c$.  Thus
$\mc A^\sharp$ is Connes-amenable by Theorem~\ref{conn_inje},
which implies in particular that $\mc A$ is unital.
Hence, if we wish to unitise $\mc A$ by using $\mc A^{cc}$,
then we need to consider ``smaller'' representations.
\noproof\end{remark}

\begin{proposition}\label{when_wap_ca}
Let $\mc A$ be a Banach algebra.  Then the following are
equivalent:
\begin{enumerate}
\item $\wap(\mc A')'$ is Connes-amenable;
\item whenever $\pi:\mc A\rightarrow\mc B(E)$ is a continuous
   representation on a reflexive Banach space $E$, there exists
   a quasi-expectation $\mc Q:\mc B(E)\rightarrow\pi(\mc A)^c$.
\end{enumerate}
\end{proposition}
\begin{proof}
Let $\pi:\mc A\rightarrow\mc B(E)$ is a continuous representation
on a reflexive Banach space $E$, so that by Proposition~\ref{wap_reps},
there is a unique weak$^*$-continuous representation $\hat\pi:
\wap(\mc A')'\rightarrow\mc B(E)$ extending $\pi$.  It is hence
sufficient to show that $\pi(\mc A)^c = \hat\pi(\wap(\mc A')')^c$.
Clearly $\pi(\mc A)^c \supseteq \hat\pi(\wap(\mc A')')^c$, while
conversely, let $T\in\pi(\mc A)^c$, so that
\[ \ip{T\pi(a)}{\tau} = \ip{\pi(a)T}{\tau}
\qquad (a\in\mc A, \tau\in E'\proten E). \]
Then let $\Phi\in\wap(\mc A')'$ and let $(a_\alpha)$ be a bounded
net in $\mc A$ which converges to $\Phi$ in the weak$^*$-topology
on $\wap(\mc A')'$.  Then, for $x\in E, \mu\in E'$ and
$T\in\pi(\mc A)^c$,
\begin{align*} &\ip{\mu}{T\hat\pi(\Phi)(x)}
= \ip{\Phi}{\pi_*( T'(\mu) \otimes x) }
= \lim_\alpha \ip{\pi_*( T'(\mu) \otimes x) }{a_\alpha} \\
&\hspace{2ex}= \lim_\alpha \ip{\mu}{T\pi(a_\alpha)(x)}
= \lim_\alpha \ip{\mu}{\pi(a_\alpha)T(x)} \\
&\hspace{2ex}= \lim_\alpha \ip{\pi_*( \mu \otimes T(x)) }{a_\alpha}
= \ip{\Phi}{\pi_*( \mu \otimes T(x)) }
= \ip{\mu}{\hat\pi(\Phi)T(x)}, \end{align*}
so that $T\in\hat\pi(\wap(\mc A')')^c$, as required.
\end{proof}

\section{WAP-compactifications for semigroups}

Semigroup algebras fit very nicely into our framework,
and the theory is well-explored.  Here we shall sketch some
results on compactifications; for further details, see \cite{BJM}.

Let $S$ be a semigroup which is also a topological space.
Then $S$ is a \emph{semitopological semigroup} when the left and right
actions of $S$ are continuous, while $S$ is a 
\emph{topological semigroup} when the multiplication map
$S\times S\rightarrow S$ is continuous.

We write $\ell^\infty(S)$ for the commutative C$^*$-algebra
of all bounded functions on $S$.  For $s\in S$, define
$\rho_s:S\rightarrow S$ by $\rho_s(t) = ts$ for $t\in S$.
Define $R_s:\ell^\infty(S)\rightarrow \ell^\infty(S)$ by
\[ R_s(f) = f\circ\rho_s \qquad (s\in S, f\in\ell^\infty(S)), \]

Let $C(S) \subseteq \ell^\infty(S)$ be the space of continuous, bounded
functions on $S$.  For $f\in\ell^\infty(S)$, we say that $f$ is
\emph{weakly almost periodic}, denoted by $f\in\mc{WAP}(S)$, when $f\in C(S)$,
and $R_S(f) := \{ R_s(f) : s\in S \}$ is relatively weakly-compact in
$\ell^\infty(S)$.  As noted in \cite[Chapter~4]{BJM}, if $S_d$ denotes the
semigroup $S$ with the discrete topology, then $\mc{WAP}(S) =
\mc{WAP}(S_d) \cap C(S)$.

\begin{theorem}\label{rep_lim_con}
Let $S$ be a semitopological semigroup, and let $f\in C(S)$.
Then $f\in\mc{WAP}(S)$ if and only if
\[ \lim_{n\rightarrow\infty} \lim_{m\rightarrow\infty} f(s_mt_n) =
\lim_{m\rightarrow\infty} \lim_{n\rightarrow\infty} f(s_mt_n), \]
whenever $(s_m)$ and $(t_n)$ are sequences of distinct elements of $S$,
and the iterated limits exist.  Then $\mc{WAP}(S)$ is a translation invariant
sub-C$^*$-algebras of $\ell^\infty(S)$ which contains the constant functions.
\end{theorem}

We shall now concentrate on the case when $S$ is discrete.  We turn
$\ell^1(S)$ into a Banach algebra with the convolution product in
the usual way (see \cite{DL} for further details about such algebras).
Using the double limit criterion above, it is a simple matter to check that
$\mc{WAP}(S) = \wap(\ell^\infty(S))$ where we treat $\ell^\infty(S)$
at an $\ell^1(S)$-bimodule.  As such, we see that the Arens products
drop to a well-defined product on $\mc{WAP}(S)'$ turning $\mc{WAP}(S)'$
into a dual Banach algebra.

Similar conclusions can be drawn when $S$ is a locally compact group
(see \cite[Chapter~7]{DL} for example) but the arguments involved
are more intricate.  For example, in \cite{Ulger}, it is shown that
$\wap(L^\infty(G)) = \mc{WAP}(G)$ for a locally compact group $G$.
The argument there seems to rely upon certain properties of groups,
and it is far from clear that an analogous result will hold in situations
where $L^1(T)$ makes sense for a topological semigroup $T$.
See \cite{LL} for recent progress in the study of when $L^1(T)$
makes sense for such a semigroup $T$.

Let $S^{\mc{WAP}}$ be the character space of the C$^*$-algebra
$\mc{WAP}(S)$.  Define a map $\epsilon:S\rightarrow S^{\mc{WAP}}$ by
letting $\epsilon(s)$ be point evaluation at $s\in S$.
We may check that the product on $\mc{WAP}(S)'$ restricts to the
character space $S^{\mc{WAP}}$, so that $S^{\mc{WAP}}$ becomes a semigroup,
and $\epsilon$ becomes a homomorphism.

For a semitopological semigroup $S$, a \emph{semitopological semigroup
compactification} of $S$ is a pair $(\psi,T)$ where $T$ is a compact, Hausdorff,
semitopological semigroup, and $\psi:S\rightarrow T$ is a continuous
homomorphism with dense range.
We do not require that $\psi$ be injective, so this differs from the
notion of a compactification of a topological group.
A semitopological semigroup compactification $(\psi,T)$ if
\emph{universal} if whenever $R$ is another compact semitopological semigroup
and $\phi:S\rightarrow R$ is a homomorphism, $\phi$ factors through
$(\psi,T)$.  Clearly any two universal compactifications are isomorphic.

\begin{theorem}
With notation as above, we have that $(\epsilon,S^{\mc{WAP}})$ is a universal
semitopological semigroup compactification of $S$.
\end{theorem}

We shall henceforth drop the $\epsilon$ and write $S^{\mc{WAP}}$ for the
\emph{$\mc{WAP}$-compactification of S}.  Notice that $\mc{WAP}(S)$ is isomorphic
to $C(S^{\mc{WAP}})$ and so $\mc{WAP}(S)'$ is isomorphic to $M(S^{\mc{WAP}})$.
We may check that the product on $\mc{WAP}(S)'$ agrees with the
natural convolution product on $M(S^{\mc{WAP}})$.

\begin{example}
Let $S$ be a discrete semigroup.  Then $\ell^1(S)$ is a Banach algebra,
and $c_0(S)$ is a predual for the Banach space $\ell^1(S)$.  We have that
the product is weak$^*$-continuous if and only if $S$ is a
\emph{weakly cancellative} discrete semigroup, that is, the
left and right actions are finite-to-one maps.  This follows by
an easy calculation: see \cite[Proposition~5.1]{Connespp} for example.
\end{example}

\begin{example}
Let $S = (\mathbb N,\max)$, so that $\ell^1(S)$ is a dual Banach
algebra with predual $c_0(S)$.
Then $S^{\mc{WAP}}$ is a compact semitopological semigroup containing
$S$ as a dense subsemigroup.  We may check that $S^{\mc{WAP}}$ is equal
to $S$ with an adjoined \emph{zero}, denoted by $\infty$, which satisfies
$\infty  n = n  \infty = \infty$ for $n\in S^{\mc{WAP}}$.  The
topology is then simply the one-point compactification.
\end{example}

\begin{example}
Let $S = (\mathbb N,\min)$, so that $\ell^1(S)$ is a Banach algebra,
but as $S$ is not weakly-cancellative, $\ell^1(S)$ is not a dual
Banach algebra with respect to $c_0(S)$.  We may check that $S^{\mc{WAP}}$
is equal to $S$ with an adjoined identity, denoted again by $\infty$, so
that $\infty  n = n  \infty = n$ for $n\in S^{\mc{WAP}}$.  The
topology is again the one-point compactification.
\end{example}

For example, the $\mc{WAP}$-compactification of $(\mathbb Z,+)$ is
a much more mysterious object.

\subsection{Injectivity and semigroup algebras}

We now apply the idea of injectivity to some semigroup algebras.
As will be seen, the results we get are rather simple, while the
necessary Banach space machinery is fairly involved, all suggesting
that we really need some further tools to make this approach worthwhile.

Let $E$ be a Banach space with a normalised basis $(e_n)$.
See \cite{LT} for further details on bases in Banach spaces.
For each $n$, there is a linear functional $e_n^*\in E'$ given by
$\ip{e_n^*}{\sum_i x_ie_i} = x_n$.  By a standard renorming of $E$,
we may suppose that the projection onto the linear span of
$(e_i)_{i=1}^n$ is norm-decreasing.

Let $S=(\mathbb N,\min)$.  There is then a natural representation $\pi:
\ell^1(S)\rightarrow\mc B(E)$ given by
\[ \pi(\delta_n)(e_m) = \begin{cases}
e_m & : m\leq n, \\ 0 &: m>n, \end{cases}
\qquad (n,m\in\mathbb N) \]
and linearity.  That is, $\pi(\delta_n)$ is the projection onto
the linear span of the first $n$ basis elements.  Let
$\mc A = \pi(\ell^1(S)) \subseteq \mc B(E)$.

Each element of $\mc B(E)$ has a natural representation as a matrix
with respect to the basis $(e_n)$.  We claim that $\mc A^c$ is just
the diagonal matricies in $\mc B(E)$.  Clearly a diagonal matrix
is in $\mc A^c$, while conversely, as $\pi(n)-\pi(n-1) = e_n^*\otimes e_n$,
we see that for $T\in\mc A^c$,
\[ T (e_n^*\otimes e_n)(e_m) = \delta_{n,m} T(e_n)
= (e_n^*\otimes e_n) T(e_m) = \ip{e_n^*}{T(e_m)} e_n
\qquad (n,m\in\mathbb N), \]
where $\delta$ here denotes the Kronecker delta.  We hence see
that $T(e_n) \in \mathbb C e_n$, for each $n\in\mathbb N$, as required.

We now claim that if there exists a quasi-expectation $\mc Q:
\mc B(E)\rightarrow\mc A^c$, then $\mc Q$ must be the canonical
projection onto the diagonal of $\mc B(E)$.  Let $n,m\in\mathbb N$,
let $a=e_n^*\otimes e_n, b=e_m^*\otimes e_m \in \mc A^c$ and let
$T\in\mc B(E)$.  Then
\[ \ip{e_n^*}{T(e_m)} \mc Q(e_m^*\otimes e_n) = \mc Q(aTb)
= a \mc Q(T) b = \ip{e_n^*}{\mc Q(T)(e_m)} e_m^*\otimes e_n, \]
so that, if $n\not=m$, then $e_m^*\otimes e_n\not\in\mc A^c$, so
$\ip{e_n^*}{\mc Q(T)(e_m)} = 0$.  Thus $\mc Q(T)(e_n) \in \mathbb C e_n$
for each $n$, and we see that
\[ \ip{e_n^*}{T(e_n)} e_n^*\otimes e_n =
\ip{e_n^*}{T(e_n)} \mc Q(e_n^*\otimes e_n)
= \ip{e_n^*}{\mc Q(T)(e_n)} e_n^*\otimes e_n, \]
and hence we see that $\mc Q(T)(e_n) = \ip{e_n^*}{T(e_n)}$, as
required.

\begin{theorem}
Let $S=(\mathbb N,\min)$.  Then $\wap(\ell^1(S)')'$ is not Connes-amenable.
\end{theorem}
\begin{proof}
By Proposition~\ref{when_wap_ca} (but really by Proposition~\ref{ca_inj})
it suffices to find a reflexive Banach space with a basis $(e_n)$ such that
the canonical projection from $\mc B(E)$ onto its diagonal is not bounded.

Following \cite[Proposition~2.b.11]{LT}, there exists a sequence
$(\beta_n)_{n\in\mathbb N}$ of positive reals tending to infinity
such that, for each $n\in\mathbb N$,
we can find a Banach space $H_n$ such that:
\begin{enumerate}
\item there is an isomorphism $\phi_n:H_n\rightarrow \ell^2_{2n}$
   (that is, the space $\mathbb C^{2n}$ with the usual Euclidean norm)
   such that $\|\phi_n\| \|\phi_n^{-1}\| \leq K$ for some absolute
   constant $K>0$;
\item $H_n$ has a normalised basis (in our sense, as above) $(e_k)_{k=1}^{2n}$;
\item there exists $(a_k)_{k=1}^{2n} \subseteq \mathbb C$ with
   \[ \Big\| \sum_{k=1}^n a_{2k-1} e_{2k-1} \Big\| \geq \beta_n
   \Big\| \sum_{k=1}^{2n} a_k e_k \Big\|. \]
\end{enumerate}
Let $K_n$ be the subspace spanned by $(e_{2k-1})_{k=1}^n$, so that
as $H_n$ is isomorphic to $\ell^2_{2n}$, let $P_n$ be the orthogonal projection
onto $K_n$ (pulled back by $\phi_n$), so that $\|P_n\| \leq K$.
Then $P_n(e_k) = e_k$ when $k$ is odd, while for even $k$, clearly
$\ip{e_k^*}{P_n(e_k)}=0$, as $P_n(e_k) \in K_n$.  Let $\mc Q_n$ be
the canonical projection of $\mc B(H_n)$ onto its diagonal, so that
$\mc Q_n(P_n) = \sum_{k=1}^n e_{2k-1}^* \otimes e_{2k-1}$.  Then
let $x = \sum_{k=1}^{2n} a_k e_k$, so that
\[ \big\| \mc Q_n(P_n)(x) \big\| = \Big\| \sum_{k=1}^n a_{2k-1} e_{2k-1}
\Big\| \geq \beta_n \|x\|, \]
so $\| \mc Q_n(P_n) \| \geq \beta_n$, and hence $\| \mc Q_n \| \geq
K^{-1} \beta_n$.

Let $E$ be the $\ell^2$-direct-sum of the $H_n$, so that if $\mc Q$ is
the canonical projection from $\mc B(E)$ onto its diagonal, we see that
$\| \mc Q \|\geq\beta_n$ for every $n$, which gives a contradiction.
We hence see that $\wap(\ell^1(S)')'$ is not Connes-amenable when
$S = (\mathbb N,\min)$.
\end{proof}

Notice that when $\pi:\ell^1(\mathbb N,\min)\rightarrow
\mc B(E)$ is a representation, then $\pi$ is a well-defined linear
operator $\ell^1(\mathbb N,\max)\rightarrow\mc B(E)$, and it is easily
checked that the map $\ell^1(\mathbb N,\max)\rightarrow \mc B(E); 
\delta_n\mapsto I_E-\pi(n)$ is a homomorphism.  The commutant for either
of these maps is equal, and hence we see that
$\wap(\ell^1(\mathbb N,\max)')'$ is also not Connes-amenable.

It seems quite possible that various interesting Banach spaces
will be generated by starting with a complex semigroup $S$ for which
we know that $\wap(\ell^1(S)')'$ is not Connes-amenable, and then
looking at the representations generated by our results.

We now briefly mention how to use the notion of injectivity to show
that a dual Banach algebra is Connes-amenable.  For example, let
$\pi:\ell^1(\mathbb Z)\rightarrow\mc B(E)$ be a representation on some
reflexive Banach space $E$.  We construct a quasi-expectation
$\mc Q:\mc B(E)\rightarrow \pi(\ell^1(\mathbb Z))^c$ by
\[ \mc Q(T) = \lim_{n\rightarrow\infty} \frac1n
\sum_{k=1}^n \pi(\delta_{-k}) T \pi(\delta_k), \]
where the limit is in the weak$^*$-topology on $\mc B(E)$.  Then,
for $m\geq 0$,
\begin{align*}
\big\| \pi(\delta_m)\mc Q(T) &- \mc Q(T)\pi(\delta_m) \big\| \\
&= \lim_{n\rightarrow\infty} \Big\| \frac1n
   \sum_{k=1}^n \big( \pi(\delta_{m-k}) T \pi(\delta_k)
   - \pi(\delta_{-k})T\pi(\delta_{k+m}) \big) \Big\| \\
&= \lim_{n\rightarrow\infty} \Big\| \frac1n
   \sum_{k=1}^m \pi(\delta_{m-k}) T \pi(\delta_k) \Big\|
\leq \lim_{n\rightarrow\infty} \frac{m}{n} \|T\| \|\pi\|^2 = 0,
\end{align*}
and a similar argument holds when $m<0$.  Thus $\mc Q$ is a projection
onto $\pi(\ell^1(\mathbb Z))^c$, and it is simple to verify that $\mc Q$
is a quasi-projection.  Of course, in this argument, we have really
used, rather directly, the fact that $\mathbb Z$ is an amenable group.

It seems more natural and profitable to study the Connes-amenability
of algebras of operators via injectivity, something we hope to
pursue in future research.

\section{Tensor products}

In this section, we shall sketch some ideas about tensor products
of Banach algebras which behave well with respect to weakly almost
periodic functionals, and then go on to give a theory of tensor
products of dual Banach algebras.

We start by sketching some results on tensor products of Banach
spaces.  We follow the notation used in Banach space theory,
namely that found in \cite{Ryan}, and in \cite{DF} (except where
this clashes with notation in \cite{Ryan}).  Note that this notation
is different to that found in \cite{tak}, for example.

Let $E$ and $F$ be Banach spaces.  We have previously defined the
projective tensor product $E\proten F$.  The \emph{injective tensor
norm}, $\epsilon$, is defined by, for $\tau = \sum_{k=1}^n x_k\otimes y_k
\in E\otimes F$,
\[ \epsilon(\tau) = \sup\Big\{ \Big| \sum_{k=1}^n
\ip{\mu}{x_k}\ip{\lambda}{y_k} \Big| : \mu\in E', \lambda\in F',
\|\mu\| = \|\lambda\|=1 \Big\}. \]
We write $E\inten F$ for the completion of $E\otimes F$ with respect
to $\epsilon$.  Notice that if we identify $E\otimes F$ with a subspace
of the finite rank operators $\mc F(E',F)$, then $\epsilon$ is the
norm induced by the operator norm under this identification.

Let $\alpha$ be a norm on $E\otimes F$.  Then $\alpha$ is a
\emph{reasonable crossnorm} when $\epsilon(\tau) \leq \alpha(\tau)
\leq \pi(\tau)$ for each $\tau\in E\otimes F$.  In this case, clearly
$\alpha(x\otimes y)=\|x\| \|y\|$ for $x\in E$ and $y\in F$.
We write $E\proten_\alpha F$ for the completion of $E\otimes F$ with
respect to $\alpha$.

A \emph{uniform crossnorm} is an assignment, to each pair of Banach
spaces $E$ and $F$, of a norm $\alpha$ with the following mapping
property.  Let $G$ and $H$ be Banach spaces, let $T\in\mc B(E,G)$ and
$S\in\mc B(F,H)$, and define $T\otimes S:E\otimes F\rightarrow G\otimes
H$ by $(T\otimes S)(x\otimes y) = T(x)\otimes S(y)$ and linearity.
Then we insist that $T\otimes S$ extends by continuity to a bounded
linear map $E\proten_\alpha F\rightarrow G\proten_\alpha H$ with
norm $\|T\| \|S\|$.  In the special case when this mapping property
holds with $E=G, F=H$ and $\alpha$ a reasonable crossnorm on $E\otimes F$,
we say (in a non-standard way) that $\alpha$ is a \emph{quasi-uniform crossnorm}.

We define an action of $E'\otimes F'$ on $E\otimes F$ by setting
\[ \ip{\mu\otimes\lambda}{x\otimes y} = \ip{\mu}{x}\ip{\lambda}{y}
\qquad (x\in E, y\in F, \mu\in E',\lambda\in F'), \]
and extending by linearity.
We define the dual norm $\alpha^s$ on $E'\otimes F'$ by
\[ \alpha^s(\sigma) = \sup\{ |\ip{\sigma}{\tau}| : \tau\in E\otimes F,
\alpha(\tau)\leq 1 \} \qquad (\sigma\in E'\otimes F'). \]
Then it may be checked that $\alpha^s$ is a reasonable crossnorm when
$\alpha$ is, and similarly for uniform crossnorms.  When $E$ and $F$
are reflexive, and $\alpha$ is a quasi-uniform crossnorm, then $\alpha^s$
is also quasi-uniform.  Then $\pi^s = \epsilon$ for all Banach spaces,
but $\epsilon^s = \pi$ only in special cases.

A \emph{tensor norm} is then a uniform crossnorm which respects
finite-dimensional subspaces in a certain sense.  We shall not
have use of this idea, but do note that many of the norms we
construct in this section are not as well behaved as those studied
in \cite{Ryan} and \cite{DF}.

As explained before, it is standard that $(E\proten F)' = \mc B(E,F')$
for Banach spaces $E$ and $F$.  For any reasonable crossnorm $\alpha$
on $E\otimes F$, as the formal inclusion map $E\proten F\rightarrow
E\proten_\alpha F$ is norm-decreasing, we may use the adjoint to
identity $(E\proten_\alpha F)'$ with a subspace of $\mc B(E,F')$,
together with the dual norm.  For example, the dual of $E\inten F$
is $\mc I(E,F')$ the \emph{integral operators} from $E$ to $F'$
(see \cite[Chapter~3]{Ryan} for further details).

As explained above, $\epsilon^s=\pi$ only in special cases, which
means that in general, the natural map $E'\proten F'\rightarrow
(E\inten F)' = \mc I(E,F')$ is only norm-decreasing.  An important
special case is when $E$ or $F$ has the \emph{metric approximation
property}, in which case $E'\proten F'$ is, isometrically, a closed
subspace of $\mc I(E,F')$.  See \cite[Section~4]{Ryan} for further
details.  Another way to state this result is to consider the natural
map from $E'\proten F'$ to $\mc B(E,F')$, which has range
$\mc N(E,F')$, the \emph{nuclear operators}.  Thus $\mc N(E,F')$
is closed in $\mc I(E,F')$ when $E$ or $F$ has the metric
approximation property.

\begin{example}
Let $X$ and $Y$ be locally compact Hausdorff spaces.  Then
$C_0(X) \inten C_0(Y) = C_0(X\times Y)$ under the obvious identification
(see \cite[Section~3.2]{Ryan}).
As $C(X)$ has the metric approximation property, we have that
$M(X) \proten M(Y)$ forms a closed subspace of $M(X\times Y)$.
We shall see below that we can fail to have equality.

Similarly, let $\mu$ and $\nu$ be measures.  Then $L^1(\mu)
\proten L^1(\nu) = L^1(\mu\times\nu)$ under the obvious
identification (see \cite[Chapter~2]{Ryan}).

Hence, if $X$ and $Y$ are discrete sets, then $c_0(X)\inten
c_0(Y) = c_0(X\times Y)$ and $c_0(X)'\proten c_0(Y)' =
\ell^1(X) \proten \ell^1(Y) = \ell^1(X\times Y) = c_0(X\times Y)'$.
\end{example}

\subsection{Tensor products of algebras}

\begin{definition}
Let $\mc A$ and $\mc B$ be Banach algebras, and define an algebra
product on $\mc A\otimes\mc B$ by $(a\otimes b)(c\otimes d) =
ac\otimes bd$, and linearity, for $a,c\in\mc A$ and $b,d\in\mc B$.
Then an \emph{algebra crossnorm} on $\mc A\otimes\mc B$ is a reasonable
crossnorm $\alpha$ such that $\mc A\proten_\alpha\mc B$ becomes
a Banach algebra.
\end{definition}

Notice that the projective tensor norm is always an algebra crossnorm,
but that the injective tensor norm may not be (indeed, it is shown in
\cite{carne} that only four of Grothendieck's fourteen ``natural''
tensor norms are always algebra crossnorms).

We have shown that WAP algebras are isomorphic (but maybe not isometric)
to closed subalgebras of $\mc B(E)$ for suitable reflexive $E$.  In this
section, it is convenient to suppose that a WAP algebra $\mc A$ is
isometric to a closed subalgebra of $\mc B(E)$ for suitable $E$.  This
can clearly be achieved by considering a suitable renorming of $\mc A$.

\begin{definition}
Let $\mc A$ and $\mc B$ be WAP algebras.  A \emph{WAP-crossnorm} on
$\mc A\otimes\mc B$ is an algebra crossnorm $\alpha$ such that, if we
form the natural chain of natural inclusion maps
$\wap(\mc A')\otimes\wap(\mc B') \subseteq \mc A'\otimes\mc B'
\subseteq \mc A'\proten_{\alpha^s}\mc B'
\subseteq (\mc A\proten_\alpha\mc B)'$, then we actually map into
$\wap((\mc A\proten_\alpha\mc B)')$, and that furthermore,
$\wap(\mc A')\otimes\wap(\mc B')$ is norming for $\mc A\proten_\alpha\mc B$.
\end{definition}

We shall see that the condition on $\wap(\mc A')\otimes\wap(\mc B')$
is natural when we come to consider dual Banach algebras.

\begin{example}\label{DBA_not_ten}
Let $K$ and $L$ be compact Hausdorff spaces, and consider the injective
tensor product $C(K)\inten C(L) = C(K\times L)$.
As these are commutative C$^*$-algebras, they are
Arens regular, and so $\wap(C(K)') = M(K)$, $\wap(C(L)') = M(L)$, and
hence $M(K)\otimes M(L) \subseteq \wap(C(K\times L))$.  As explained above,
we induce the projective tensor norm on $M(K)\otimes M(L)$ by
embedding it into $C(K\times L)' = M(K\times L)$.  Then
$(M(K)\proten M(L))' = \mc B(M(K),M(L)')$, and so we get the chain of
isometric inclusions
\[ C(K) \inten C(L) \subseteq \mc A(M(K),C(L))
\subseteq \mc A(M(K),M(L)') \subseteq \mc B(M(K),M(L)'). \]
Thus $M(K)\otimes M(L)$ is norming for $C(K\times L)$.  We could also
see this directly by considering point masses in $M(K)$ and $M(L)$.

However, $M(K)\otimes M(L)$ is not in general dense in $M(K\times L)$.
Let $L$ be a compact Hausdorff space such that $M(L)$ does not have the
\emph{Radon-Nikod\'{y}m property} (see \cite{DU} for what this technical
condition is).  For example, \cite[Chapter~VII]{DU} shows that this
holds when there is a separable subspace of $C(L)$ without a separable
dual.  As indicated in \cite[Chapter~VI, Corollary~6]{DU}, there then exists
an integral, non-nuclear operator from $C(K)$ to $M(L)$ whenever $K$
contains a \emph{perfect} subset (that is, a closed subset with no
isolated points).  For example, let $\mathbb T = \{z\in\mathbb C :
|z|=1\}$, so that $\mathbb T$ is perfect, and $M(\mathbb T)$ does not
have the Radon-Nikod\'{y}m property.  Thus
$M(\mathbb T) \otimes M(\mathbb T)$ is not dense
$(C(\mathbb T) \inten C(\mathbb T))'$.
\end{example}

\begin{example}
Let $S$ and $T$ be discrete semigroups, so that $\wap(\ell^1(S)') =
C(S^{\mc WAP})$.  Then $\ell^1(S) \proten \ell^1(T) = \ell^1(S\times T)$
as a Banach algebra, and so $\wap( \ell^1(S\times T)' )
= C((S\times T)^{\mc WAP})$.  We claim that $(S\times T)^{\mc WAP}
= S^{\mc WAP} \times T^{\wap}$, which follows easily by the universality
property of $(S\times T)^{\mc WAP}$.  We then see that
\begin{align*}
\wap(\ell^1(S)') \inten \wap(\ell^1(T)') &=
C(S^{\mc WAP}) \inten C(T^{\mc WAP}) =
C( S^{\mc WAP} \times T^{\mc WAP} ) \\
&= C( (S\times T)^{\mc WAP} ) =
\wap( (\ell^1(S)\proten\ell^1(T))' ), \end{align*}
so that certainly $\wap(\ell^1(S)') \otimes \wap(\ell^1(T)')$ is
dense in $\wap( (\ell^1(S)\proten\ell^1(T))' )$.

As before, this argument also works for general locally compact
groups $G$ and $H$.
\end{example}

\begin{example}\label{semigroups_ten}
Let $S$ and $T$ be discrete weakly cancellative semigroups,
so that $c_0(S)\inten c_0(T) = c_0(S\times T)$, and hence
$(c_0(S)\inten c_0(T))' = \ell^1(S\times T) = \ell^1(S)\proten\ell^1(T)$
is a dual Banach algebra.

Similarly, let $S$ and $T$ be locally compact groups, so that
$M(S)$ has the predual $C_0(S)$, and similarly for $M(T)$
(see \cite{Runde1}).  Then, as above, we see that $M(S)\proten
M(T)$ is a closed, norming (and hence weak$^*$-dense) subspace of
$M(S\times T)=C_0(S\times T)' = ( C_0(S)\inten C_0(T) )'$.
\end{example}

\begin{proposition}\label{when_proj_wapnorm}
Let $\mc A$ and $\mc B$ be Banach algebras.  Then the natural map
$\wap(\mc A')\otimes\wap(\mc B')\rightarrow\mc A'\proten_{\pi^s}
\mc B'$ maps into $\wap((\mc A\proten\mc B)')$.  When $\mc A$
and $\mc B$ are Arens regular and one has the metric approximation
property, $\pi$ is a WAP-crossnorm on $\mc A\otimes\mc B$.
\end{proposition}
\begin{proof}
Recall that $\pi^s=\epsilon$, and let $\mu\in\wap(\mc A')$ and
$\lambda\in\wap(\mc B')$.  Then there exists a reflexive Banach
space $E$ and maps $S_\mu:\mc A\rightarrow E$ and $T_\mu:E
\rightarrow\mc A'$ such that $T_\mu S_\mu=L_\mu$, that is, $T_\mu
S_\mu(a)=a\cdot\mu$ for $a\in\mc A$.  Similarly there exists a
reflexive Banach space $F$ and maps $S_\lambda$ and $T_\lambda$.
We shall see below that there exists a uniform-crossnorm $\alpha$
such that $E\proten_\alpha F$ is reflexive.  As $\epsilon\leq
\alpha\leq\pi$, we see that
we can factor the map $L_\mu\otimes L_\lambda:\mc A\proten\mc B
\rightarrow \mc A'\inten\mc B'$ as
\[ \xymatrix{ \mc A\proten\mc B \ar[r]^{S_\mu\otimes S_\lambda}
& E\proten F \ar[r] & E\proten_\alpha F 
\ar[r]^{T_\mu\otimes T_\lambda} & \mc A'\proten_\alpha \mc B'
\ar[r] & \mc A'\inten \mc B',} \]
and so $L_\mu\otimes L_\lambda = L_{\mu\otimes\lambda}$ is weakly-compact.
By linearity, the argument is complete.

Suppose now that $\mc A$ has the metric approximation property.
By \cite[Theorem~4.14]{Ryan}, the canonical map $\mc A\proten\mc B
\rightarrow (\mc A'\inten\mc B')'$ is an isometry, so that
$\mc A'\inten\mc B'$ is norming for $\mc A\proten\mc B$.  As $\mc A$
and $\mc B$ are Arens regular, we see that $\pi$ is indeed a
WAP-crossnorm on $\mc A\otimes\mc B$.  The case for $\mc B$
follows by symmetry.
\end{proof}

\begin{example}
Let $c = C(\mathbb N_\infty)$ be the space of convergent sequences,
where $\mathbb N_\infty$ is the one-point compactification of
$\mathbb N$.  Then $c'=\ell^1$ naturally, and so $(c\proten c)' =
\mc B(c,\ell^1) = \mc A(c,\ell^1) = \ell^1\inten\ell^1$.  Then $\pi$
is a WAP-crossnorm for $c\otimes c$ as $c$ is Arens regular.  In fact,
$c\proten c$ is Arens regular.

For general compact Hausdorff spaces $K$ and $L$, by the above
proposition, we see that $\pi$ is a WAP-crossnorm on $C(K)\otimes
C(L)$.  However, as noted in \cite{UlgerCorr}, when $G$ is a compact
group, the algebra $C(G)\proten C(G)$ contains a copy of $A(G)$, the
Fourier algebra of $G$, so that $C(G)\proten C(G)$ is not Arens regular.
Hence $M(G) \otimes M(G) \subseteq (C(G)\proten C(G))'$ cannot
be dense.
\noproof\end{example}

See \cite{Ulger2} for details about when $\mc A\proten\mc B$ is
Arens regular, but bear in mind the correction \cite{UlgerCorr}.

It seems that in general, $\wap(\mc A')\inten\wap(\mc B')$ need not
be norming for $\mc A\proten\mc B$.  Also, we see no way to adapt
the above proof to the case $\mc A\proten_\beta\mc B$ for an arbitrary
algebra-crossnorm $\beta$.

\begin{proposition}\label{ref_ten_norm}
Given reflexive Banach spaces $E$ and $F$, there exists a tensor
norm $\alpha$ on $E\otimes F$ such that $E\proten_\alpha F$ is
reflexive, and such that $E'\otimes F'$ is dense in
$(E\proten_\alpha F)'$.
\end{proposition}
\begin{proof}
For example, for $1<p<\infty$, let $g_p$ and $d_p$ be the
\emph{Chevet-Saphar} tensor norms, as defined in \cite[Chapter~6]{Ryan}
or \cite{as}.  Then, by \cite[Corollary~3.2]{as}, we have that
$E\proten_{g_p} F$ is reflexive, whenever $E$ and $F$ are, and
similarly for $d_p$.  Furthermore, $E'\otimes F'$ is indeed dense
in $(E\proten_{g_p} F)'$ and $(E\proten_{d_p} F)'$.  As $E'\otimes F'$
is dense in $(E\proten_\alpha F)'$, we see that
$(E\proten_\alpha F)' = E'\proten_{\alpha^s}F'$.
\end{proof}

\begin{remark}
One disadvantage of the Chevet-Saphar tensor norms is that they
are not symmetric, in that $E\proten_{g_p} F$ is not in general
isomorphic to $F\proten_{g_p} E$ (in contrast to the injective or
projective tensor norms).  However, $d_2$ and $g_2$ behave rather
nicely, in that $d_2^s = g_2$, $g_2^s = d_2$, and $H\proten_{d_2} K
= H\proten_{g_2} K = H\otimes_2 K$ whenever $H$ and $K$ are Hilbert
spaces (here $H\otimes_2 K$ is the usual Hilbertian tensor product).
If we interpolate between $d_2$ and $g_2$ then we can verify that we
end up with a symmetric tensor norm $\alpha$ such that $E\proten_\alpha F$ is
reflexive when $E$ and $F$ are, and such that $H\proten_\alpha K
= H\otimes_2 K$ for Hilbert spaces $H$ and $K$.
\noproof\end{remark}

Let $\mc A$ be a Banach algebra such that there is a reflexive
Banach space $E$ and an isometric representation
$\pi_{\mc A}:\mc A\rightarrow\mc B(E)$ such that
$\pi_{\mc A}'\kappa_{E'\proten E}$ takes the unit
ball of $E'\proten E$ onto the unit ball of $\wap(\mc A')$.
From our previous work, this is equivalent to $\wap(\mc A')$ being norming
for $\mc A$.  Suppose that $\mc B$ is similar, with
$\pi_{\mc B}:\mc B\rightarrow\mc B(F)$, say.
Now let $\alpha$ be some quasi-uniform crossnorm satisfying the conclusions
of Proposition~\ref{ref_ten_norm}.
Then $\pi = \pi_{\mc A}\otimes\pi_{\mc B}:\mc A\otimes
\mc B\rightarrow \mc B(E\proten_\alpha F)$, defined by
\[ \pi(a\otimes b)(x\otimes y) = \pi_{\mc A}(a)(x) \otimes \pi_{\mc B}(b)(y)
\qquad (a\in\mc A,b\in\mc B,x\in E,y\in F), \]
is a representation.  Use this to induce a norm
$\|\cdot\|_{\pi,\alpha}$ on $\mc A\otimes\mc B$, and denote the
completion by $\mc A \proten_{\pi,\alpha} \mc B$.

\begin{proposition}\label{rep_wap_crnorm}
With notation as above, $\|\cdot\|_{\pi,\alpha}$ is a WAP-crossnorm
on $\mc A\otimes\mc B$.
\end{proposition}
\begin{proof}
By assumption, $\|a\otimes b\|_{\pi,\alpha} = \|\pi_{\mc A}(a)\|
\|\pi_{\mc B}(b)\| = \|a\| \|b\|$ for $a\in\mc A$ and $b\in\mc B$, so the
triangle-inequality implies that $\|\cdot\|_{\pi,\alpha}\leq\pi(\cdot)$.
Let $\lambda_{\mc A}\in\mc A'$ and $\lambda_{\mc B}\in\mc B'$ be
such that $\|\lambda_{\mc A}\| = \|\lambda_{\mc B}\| = 1$.  As $\wap(\mc A')$
is norming for $\mc A$, the unit ball of $\wap(\mc A')$ is weak$^*$-dense
in the unit ball of $\mc A'$, and so there exists a net $(\mu^{\mc A}_\alpha)$
in $\wap(\mc A')$ such that $\|\mu^{\mc A}_\alpha\|\leq\|\lambda_{\mc A}\|$
for each $\alpha$, and $\lim_\alpha \ip{\mu^{\mc A}_\alpha}{a} =
\ip{\lambda_{\mc A}}{a}$ for $a\in\mc A$.  By the assumption on $\pi_{\mc A}$,
for each $\alpha$, we can find $\sigma^{\mc A}_\alpha\in E'\proten E$
with $\|\sigma^{\mc A}_\alpha\| = \|\mu^{\mc A}_\alpha\|$ and
$\pi_{\mc A}'\kappa_{E'\proten E}(\sigma^{\mc A}_\alpha) = 
\mu^{\mc A}_\alpha$.  Similarly, we can find
$(\sigma^{\mc B}_\beta) \subseteq F'\proten F$ for $\lambda_{\mc B}$.
Then, for $\tau=\sum_{k=1}^n a_k\otimes b_k\in \mc A\otimes\mc B$,
\begin{align*} &\ip{\lambda_{\mc A}\otimes\lambda_{\mc B}}{\tau} =
\sum_{k=1}^n \lim_\alpha \ip{\mu^{\mc A}_\alpha}{a_k}
   \lim_\beta \ip{\mu^{\mc B}_\beta}{b_k} \\
&= \sum_{k=1}^n \lim_\alpha \ip{\pi_{\mc A}(a_k)}{\sigma^{\mc A}_\alpha}
   \lim_\beta \ip{\pi_{\mc B}(b_k)}{\sigma^{\mc B}_\beta}
= \lim_\alpha \lim_\beta \ip{\pi(\tau)}{\sigma^{\mc A}_\alpha
   \otimes \sigma^{\mc B}_\beta}.
\end{align*}
By taking the supremum over $\|\lambda_{\mc A}\| = \|\lambda_{\mc B}\|
= 1$, we consequently conclude that $\epsilon(\cdot)\leq\|\cdot\|_{\pi,\alpha}$.
Thus $\|\cdot\|_{\pi,\alpha}$ is a reasonable-crossnorm, and so
clearly it is an algebra-crossnorm.

Now let $G = E\proten_\alpha F$ and $\mc C = \mc A\proten_{\pi,\alpha}\mc B$.
We may treat $\pi$ as an isometry from $\mc C$ into $\mc B(G)$, and so
$\pi'\kappa_{G'\proten G}$ is a norm-decreasing map from $G'\proten G$
to $\mc C'$.  Define
\[ X = \big\{ \pi'\kappa_{G'\proten G}((\mu_E\otimes\mu_F)\otimes(x_E\otimes x_F))
: \mu_E\in E', \mu_F\in F', x_E\in E, x_F\in F \big\}, \]
so that clearly $X \subseteq \wap( \mc C' )$.  For $\mu_E\in E', \mu_F\in F',
x_E\in E$ and $x_F\in F$, we have that for $a\in\mc A$ and $b\in\mc B$,
\[ \ip{\pi'\kappa_{G'\proten G}((\mu_E\otimes\mu_F)\otimes(x_E\otimes x_F))}{a\otimes b}
= \ip{\pi_{\mc A}'(\mu_E\otimes x_E)}{a} \ip{\pi_{\mc B}'(\mu_F\otimes x_F)}{b}. \]
Hence, by the assumptions on $\pi_{\mc A}$ and $\pi_{\mc B}$,
the linear span of $X$ is dense in $\wap(\mc A')\otimes\wap(\mc B')$,
and so $\wap(\mc A')\otimes\wap(\mc B') \subseteq \wap(\mc C')$, as required.
Finally, almost by definition, the linear span of $X$ is norming for
$\pi(\mc A\otimes\mc B)$, and so $\|\cdot\|_{\pi,\alpha}$ is a WAP-crossnorm.
\end{proof}

Conversely, suppose that $\beta$ is some WAP-crossnorm on
$\mc A\otimes\mc B$, so that there exists a reflexive Banach space $G$
and an isometric representation $\pi:\mc A\proten_\beta\mc B\rightarrow
\mc B(G)$.  However, it need not be the case that $G=E\proten_\alpha F$
for some Banach spaces $E$ and $F$, and some uniform-crossnorm $\alpha$,
with $\pi=\pi_{\mc A}\otimes\pi_{\mc B}$ for suitable $\pi_{\mc A}:
\mc A\rightarrow\mc B(E)$ and $\pi_{\mc B}:\mc A\rightarrow\mc B(F)$.

This mirrors the behaviour of C$^*$-algebras.  Recall that the
\emph{minimal} C$^*$-tensor product of two C$^*$-algebras $\mc A$ and
$\mc B$ is that defined by taking faithful $*$-representations of $\mc A$
and $\mc B$ on Hilbert spaces $H$ and $K$, respectively, and letting
$\mc A\otimes_{\min}\mc B$ be the closure of $\mc A\otimes\mc B$ in
$\mc B(E\otimes_2 F)$.  It turns out that this is independent of
the $*$-representations taken (as long as they are faithful).  The
\emph{maximal} C$^*$-tensor product is that defined by taking the
supremum over any $*$-representation of $\mc A\otimes\mc B$.  As
indicated by their names, the minimal and maximal C$^*$-tensor norms
are indeed the smallest and greatest norms on $\mc A\otimes\mc B$
which satisfy the C$^*$-condition $\|\tau^*\tau\|=\|\tau\|^2$
for $\tau\in\mc A\otimes\mc B$.  Then a 
C$^*$-algebra $\mc A$ is \emph{nuclear} if $\mc A\otimes_{\min}\mc B=
\mc A\otimes_{\max}\mc B$ for all C$^*$-algebras $\mc B$.  This is
actually equivalent to $\mc A$ being amenable (see \cite{RundeBook}).

As C$^*$-algebras are always Arens regular, we see that a C$^*$-tensor
norm $\alpha$ on $\mc A\otimes\mc B$ is a WAP-crossnorm if and only if
$\mc A'\otimes\mc B'$ is norming for $\mc A\proten_\alpha\mc B$.  However,
this is always true for $\mc A'\otimes_{\min}\mc B'$, essentially for
the same reasons as in the proof of Proposition~\ref{rep_wap_crnorm}.
By \cite[Proposition~4.10]{tak}, the norm induced on $\mc A'\otimes\mc B'$
by $\mc A\otimes_{\max}\mc B$ always agrees with that induced by
$\mc A\otimes_{\min}\mc B$, so if $\alpha$ is a C$^*$-tensor norm,
then as $\min\leq\alpha\leq\max$, we see that $\mc A'\proten_{\alpha^s}
\mc B' = \mc A'\proten_{\min^s}\mc B'$.  Hence $\alpha$ is a WAP-crossnorm,
as we might hope.

We may be tempted to define the minimal and maximal WAP-crossnorms
in a similar fashion, given that we know that at least one
WAP-crossnorm must exist.  However, it is rather unclear if such
minimal and maximal norms exist.

\begin{proposition}
Let $\mc A$ and $\mc B$ be Banach algebras.  There exists an
algebra crossnorm $\max$ on $\mc A\otimes\mc B$ such that if $\alpha$
is any WAP-crossnorm on $\mc A\otimes\mc B$, then $\alpha
\leq\max$, and such that $\mc A\proten_{\max}\mc B$ is a WAP-algebra.
\end{proposition}
\begin{proof}
Define $\max$ on $\mc A\otimes\mc B$ by $\max(\tau) = \sup\{ \|\pi(\tau)\|
\}$, where we take the supremum over all algebra homomorphisms
$\pi:\mc A\otimes\mc B \rightarrow\mc B(E)$, where $E$ is a reflexive
Banach space, and $\|\pi(a\otimes b)\| = \|a\| \|b\|$ for $a\in\mc A$
and $b\in\mc B$.  Let $\|\cdot\|_{\pi,\alpha}$ be a WAP-crossnorm given
by Proposition~\ref{rep_wap_crnorm}, so that $\epsilon \leq
\|\cdot\|_{\pi,\alpha} \leq \max$.  As $\max(a\otimes b) = \|a\| \|b\|$
for $a\in\mc A$ and $b\in\mc B$, we see that $\max\leq\pi$, so that $\max$
is a reasonable crossnorm.  Clearly then $\max$ is an algebra crossnorm.

Let $\alpha$ be some WAP-crossnorm on $\mc A\otimes\mc B$, so we
can find a representation $\pi : \mc A\proten_{\alpha}\mc B
\rightarrow\mc B(E)$ for some reflexive Banach space $E$, such that
$\pi'$ take the unit ball of $E'\proten E$ onto the unit ball of
$\wap((\mc A\proten_\alpha\mc B)')$.  As $\wap(\mc A')\otimes
\wap(\mc B') \subseteq\wap((\mc A\proten_\alpha\mc B)')$ is norming
for $\mc A\proten_\alpha\mc B$, we see that $\pi$ is an isometry.
Consequently, $\alpha\leq\max$.

For each $\tau\in\mc A\otimes\mc B$ and each $\epsilon>0$, let
$\pi_{\tau,\epsilon}:\mc A\otimes\mc B\rightarrow\mc B(E_{\tau,\epsilon})$
be some representation on a reflexive Banach space such that
$\|\pi_{\tau,\epsilon}(\tau)\| > \max(\tau)-\epsilon$.  Let $F = \ell^2
\big( \bigoplus_{\tau,\epsilon} E_{\tau,\epsilon} \big)$, so that
$\pi = \bigoplus \pi_{\tau,\epsilon}$ is a representation of $\mc A
\otimes\mc B$ on $\mc B(F)$.  Clearly then, for $\tau\in\mc A\otimes\mc B$,
$\|\pi(\tau)\| = \sup \| \pi_{\tau,\epsilon}(\tau)\| = \max(\tau)$,
and so $\pi$ extends to an isometric representation $\pi:\mc A\proten_{\max}
\mc B\rightarrow\mc B(F)$.
\end{proof}

In general, however, we see no way to show that this norm is a
WAP-crossnorm.

\subsection{Dual Banach algebras}

Recall that for two W$^*$-algebras $(\mc A,\mc A_*)$ and
$(\mc B,\mc B_*)$, we let $X = \mc A\otimes_{\min}\mc B$,
and then regard $\mc A_*\otimes\mc B_*$ as a subspace of $X'$,
with closure $\mc A_*\overline\otimes\mc B_*$.
Then $\mc A\overline\otimes\mc B$, the \emph{W$^*$-tensor product}
of $\mc A$ and $\mc B$, is the dual of $\mc A_*\overline\otimes\mc B_*$.
It may be checked that $\mc A\otimes\mc B$
becomes a weak$^*$-dense subalgebra of $\mc A\overline\otimes\mc B$.
For example, when $H$ and $K$ are Hilbert spaces, we have that
$\mc B(H)\overline\otimes\mc B(K) = \mc B(H\otimes_2 K)$; in particular,
$\mc B(H)\otimes\mc B(K)$ need not be norm-dense in 
$\mc B(H)\overline\otimes\mc B(K)$ (see \cite[Exercise~11.5.7]{kad}).

Let $(\mc A,\mc A_*)$ and $(\mc B,\mc B_*)$ be dual Banach algebras,
and let $\beta$ be a WAP-crossnorm on $\mc A\otimes\mc B$.  By
assumption, $\wap(\mc A') \proten_{\beta^s} \wap(\mc B')$ is norming
for $\mc A\proten_\beta\mc B$ and maps into
$\wap((\mc A\proten_\beta\mc B)')$.  Hence we see that
$X = \mc A_*\proten_{\beta^s}\mc B_*$ is a closed
$\mc A\proten_\beta\mc B$-submodule of
$\wap((\mc A\proten_\beta\mc B)')$.  As such, by Proposition~\ref{wap_dba},
$(X',X)$ becomes a dual Banach algebra.  Notice that the norm on
$X=\mc A_*\proten_{\beta^s}\mc B_*$ is given by
\[ \|\tau\|_X = \sup\{ |\ip{\tau}{u}| : u\in\mc A\otimes\mc B,
\beta(u)\leq 1\} \qquad (\tau\in X). \]
It is hence immediate that $\mc A\otimes\mc B$ is weak$^*$-dense in $X'$,
so we may regard $X'$ as a dual Banach algebra
tensor product of $\mc A$ and $\mc B$, and denote $X'$ by
$\mc A\overline\otimes_\beta\mc B$.
When $\mc A_*\proten_{\beta^s}\mc B_*$ is norming for
$\mc A\proten_\beta\mc B$, we see that $\mc A\proten_\beta\mc B$
is even a closed subspace of $\mc A\overline\otimes_\beta\mc B$.

\begin{example}
Let $G$ and $H$ be locally compact groups, so that $M(G)$ is
a dual Banach algebra with predual $C_0(G)$, and similarly for
$H$ (see Example~\ref{semigroups_ten}).
By Proposition~\ref{when_proj_wapnorm}, $\wap(M(G)')\otimes
\wap(M(H)')$ maps into $\wap( (M(G)\proten M(H))' )$.  As the
map $M(G)\proten M(H) \rightarrow (C_0(G)\inten C_0(H))'$
is an isometry, we see that
\[ C_0(G) \otimes C_0(H) \subseteq C_0(G) \inten C_0(G)
\subseteq (M(G)\proten M(H))' = \mc B(M(G),M(H)') \]
is norming for $M(G)\proten M(H)$.
Thus $\pi$ is a WAP-crossnorm on $M(G)\otimes M(H)$.

We may hence form $M(G) \overline\otimes_\pi M(H)$.  By definition,
it is the weak$^*$-closure of $M(G)\otimes M(H)$ in the
dual of $C_0(G)\inten C_0(H)$, that is, in $M(G\times H)$.
As explained in Example~\ref{semigroups_ten}, $M(G)\otimes M(H)$
is certainly weak$^*$-dense in $M(G\times H)$, and so
$M(G) \overline\otimes_\pi M(H) = M(G\times H)$ as we might hope.
\end{example}

Given dual Banach algebras $(\mc A,\mc A_*)$ and $(\mc B,\mc B_*)$,
we can find weak$^*$-weak$^*$-continuous representations $\pi_{\mc A}:\mc A
\rightarrow\mc B(E)$ and $\pi_{\mc B}:\mc B\rightarrow\mc B(F)$.
We may hence define the WAP-crossnorm $\beta:=\|\cdot\|_{\pi,\alpha}$,
as in Proposition~\ref{rep_wap_crnorm}, leading to
$\mc A\proten_{\pi,\alpha}\mc B$ and hence
$\mc A\overline\otimes_\beta\mc B$.
Alternatively, we may simply define $\mc A\overline\otimes_{\pi,\alpha}\mc B$
to be the weak$^*$-closure of $\pi_{\mc A}(\mc A)\otimes\pi_{\mc B}(\mc B)$ in
$\mc B(E\proten_\alpha F)$.

\begin{proposition}\label{dba_tensor}
With notation as above, there is a natural norm-decreasing map from
$\mc A_* \proten_{\beta^s} \mc B_*$ to
$(\mc A\overline\otimes_{\pi,\alpha}\mc B)_*$ which has dense range.
If $\mc A_*\proten_{\beta^s}\mc B_*$ is norming for
$\mc A\proten_\beta\mc B$, then this map is an isometry.
\end{proposition}
\begin{proof}
We identify $\mc A$ with its weak$^*$-closed image under
$\pi_{\mc A}:\mc A\rightarrow\mc B(E)$, and hence we identify
$\mc A_*$ with a quotient of $E'\proten E$, namely
$E'\proten E / {^\perp}\mc A$, and similarly for $\mc B$.
Then $\mc A\proten_\beta\mc B$ is the closure of $\mc A\otimes\mc B$
in $\mc B(E\proten_\alpha F)$.  We have that
$(\mc A\overline\otimes_{\pi,\alpha}\mc B)_*
= \mc B(E\proten_\alpha F)_* / {^\perp}(\mc A\otimes\mc B)$
where $\mc B(E\proten_\alpha F)_* = (E'\proten_{\alpha^s}F') \proten
(E\proten_\alpha F)$.

We define a natural map $\theta:\mc A_*\otimes\mc B_* \rightarrow
(\mc A\overline\otimes_{\pi,\alpha}\mc B)_*$
as follows.  For $\tau\in\mc A_*$, $\sigma\in\mc B_*$ and
$\epsilon>0$, we have representations
\[ \tau = \sum_{n=1}^\infty x^E_n\otimes\mu^E_n + {^\perp}\mc A
\in E'\proten E / {^\perp}\mc A,
\quad \sigma = \sum_{n=1}^\infty x^F_n\otimes\mu^F_n + {^\perp}\mc B
\in F'\proten F / {^\perp}\mc B, \]
where $\sum_{n=1}^\infty \|x^E_n\|\|\mu^E_n\|<\|\tau\|+\epsilon$ and
$\sum_{n=1}^\infty \|x^F_n\|\|\mu^F_n\|<\|\sigma\|+\epsilon$.  Then let
\[ \theta(\tau\otimes\sigma) = u = \sum_{n=1}^\infty \sum_{m=1}^\infty
(x^E_n\otimes x^F_m) \otimes (\mu^E_n\otimes\mu^F_m)
+ {^\perp}(\mc A\otimes\mc B) \in(\mc A\overline\otimes_{\pi,\alpha}\mc B)_*, \]
so that $\|u\| \leq \sum_{n=1}^\infty \sum_{m=1}^\infty
\|x^E_n\| \|x^F_m\| \|\mu^E_n\| \|\mu^F_m\|
< (\|\sigma\|+\epsilon)(\|\tau\|+\epsilon)$.
For $a\in\mc A$ and $b\in\mc B$, we see that
\[ \ip{a\otimes b}{u} = \sum_{n=1}^\infty \sum_{m=1}^\infty
\ip{a}{x^E_n\otimes\mu^E_n} \ip{b}{x^F_m\otimes\mu^F_m}
= \ip{a}{\tau}\ip{b}{\sigma}
= \ip{a\otimes b}{\tau\otimes\sigma}. \]
In particular, $u$ does not depend upon the choice of representatives
for $\tau$ and $\sigma$, and so $\theta$ is well-defined.  Notice also that
for $\psi\in\mc A_*\otimes\mc B_*$, 
\begin{align*}
\beta^s(\psi) &= \sup\{ |\ip{\psi}{v}| : v\in\mc A\otimes\mc B,
\beta(v)\leq 1 \} \\&= \sup\{ |\ip{\theta(\psi)}{v}| : v\in\mc A\otimes\mc B
\subseteq\mc B(E\proten_\alpha F), \|v\|\leq 1 \} \leq \|\theta(\psi)\|,
\end{align*}
so that $\theta$ extends to a norm-decreasing map $\mc A_*\proten_{\beta^s}\mc B_*
\rightarrow(\mc A\overline\otimes_{\pi,\alpha}\mc B)_*$.
Notice that if we knew that the unit ball of
$\mc A\otimes\mc B$ were weak$^*$-dense in the unit ball of
$\mc A\overline\otimes_{\pi,\alpha}\mc B$ (that is, $\mc A_*\proten_{\beta^s}\mc B_*$
were norming for $\mc A\proten_\beta\mc B$), then we
would even have that $\theta$ was an isometry.  This is the case for
von Neumann algebras, for example, by the Kaplansky Density theorem
(\cite[Section~II, Theorem~4.8]{tak}).

We shall now show that $\mc A_*\otimes\mc B_*$ is dense in
$(\mc A\overline\otimes_{\pi,\alpha}\mc B)_*$,
which will complete the proof.
Let $\tau\in (\mc A\overline\otimes_{\pi,\alpha}\mc B)_*$, and pick a representation
\[ \tau = \sum_{n=1}^\infty u_n \otimes v_n + {^\perp}(\mc A\otimes\mc B), \]
with $(u_n)\subseteq(E'\proten_{\alpha^s}F')$ and
$(v_n)\subseteq(E\proten_\alpha F)$ satisfying $\sum_{n=1}^\infty
\alpha^s(u_n)\alpha(v_n) <\infty$.  By approximation, we may actually suppose
that $(u_n)\subseteq E'\otimes F'$ and $(v_n)\subseteq E\otimes F$.
We can then find representations
\[ u_n = \sum_{k=1}^\infty \phi^{(n)}_k \otimes \psi^{(n)}_k,
\quad v_n = \sum_{k=1}^\infty x^{(n)}_k \otimes y^{(n)}_k, \]
where for each $n$, eventually $\phi^{(n)}_k=0$, and so forth.
For each $n$, define
\[ \mu_n = \sum_{k=1}^\infty \phi^{(n)}_k \otimes x^{(n)}_k
+ {^\perp}\mc A \in \mc A_*,
\quad \lambda_n = \sum_{k=1}^\infty \psi^{(n)}_k \otimes y^{(n)}_k
+ {^\perp}\mc B \in \mc B_*, \]
noticing that each of these is a finite sum.  For $a\in\mc A$
and $b\in\mc B$, we see that
\[ \ip{a\otimes b}{\mu_n\otimes\lambda_n}
= \sum_{k=1}^\infty \ip{a}{\phi^{(n)}_k \otimes x^{(n)}_k}
\ip{b}{\psi^{(n)}_k \otimes y^{(n)}_k} =
\ip{a\otimes b}{u_n\otimes v_n}, \]
for $n\geq1$.
Consequently, for $c\in\mc A\overline\otimes_{\pi,\alpha}\mc B$ and $N\geq1$,
as $\mc A\otimes\mc B$ is weak$^*$-dense in $\mc A\overline\otimes_{\pi,\alpha}\mc B$,
\[ \Big| \ip{c}{\tau} - \sum_{n=1}^N
   \ip{c}{\mu_n\otimes\lambda_n} \Big|
= \Big| \sum_{n=N+1}^\infty \ip{c}{u_n\otimes v_n} \Big|
\leq \|c\| \sum_{n=N+1}^\infty \alpha^s(u_n)\alpha(v_n). \]
Thus $\mc A_*\otimes\mc B_*$ is indeed dense in
$(\mc A\overline\otimes_{\pi,\alpha}\mc B)_*$.
\end{proof}

As indicated, the lack of a generalisation of the Kaplansky
Density Theorem shows that in general $\mc A\overline\otimes_\beta\mc B$
and $\mc A\overline\otimes_{\pi,\alpha}\mc B$ are different.
The following provides an example of a general Banach algebra in
which the theory works well.

\begin{proposition}\label{for_op_spaces}
Let $E$ and $F$ be reflexive Banach spaces, and let $\alpha$
be a tensor norm on $E\otimes F$ satisfying the conclusions of
Proposition~\ref{ref_ten_norm}.  Form the tensor products
$\mc B(E)\overline\otimes_\beta\mc B(F)$ and $\mc B(E)\overline\otimes_{
\pi,\alpha}\mc B(F)$ by using the trivial representations of
$\mc B(E)$ on itself, and the same for $\mc B(F)$.
Then these tensor products agree with $\mc B(E\proten_\alpha F)$.
\end{proposition}
\begin{proof}
By the above proposition, $\mc B(E)_* \otimes \mc B(F)_*$ is dense in
the predual of $\mc B(E)\overline\otimes_{\pi,\alpha}\mc B(F)$.
Furthermore, $E'\otimes E\otimes F'\otimes F$ is dense in both
$\mc B(E)_* \otimes \mc B(F)_*$ and $\mc B(E\proten_\alpha F)_*$.
So if our natural map is an isometry, that is, $\mc B(E)_*\proten_{\beta^s}
\mc B(F)_*$ is norming for $\mc B(E)\proten_\beta\mc B(F)$,
then the proof is complete.

If we identify $\mc B(E)\otimes\mc B(F)$ as a subalgebra of
$\mc B(E\proten_\alpha F)$ then $\beta$ agrees with the
operator norm.
Let $u\in\mc B(E)\otimes\mc B(F)$ and $\epsilon>0$,
so we may find $\sigma\in E'\otimes F'$ and $\tau\in E\otimes F$
with $\alpha^s(\sigma)\leq1$, $\alpha(\tau)\leq 1$ and
$|\ip{\sigma}{u(\tau)}|>\beta(u)-\epsilon$.  Let
\[ \sigma=\sum_{i=1}^n \mu_i\otimes\lambda_i,\quad
\tau=\sum_{j=1}^m x_j\otimes y_j, \]
and define
\[ v = \sum_{i=1}^n \sum_{j=1}^m (\mu_i\otimes x_j) \otimes
(\lambda_i\otimes y_j) \in (E'\otimes E)\otimes (F'\otimes F)
\subseteq \mc B(E)_* \otimes \mc B(F)_*. \]
A simple calculation shows that $\ip{w}{v} = \ip{\sigma}{w(\tau)}$
for any $w\in\mc B(E)\otimes\mc B(F)$, so that
\begin{align*}
\beta^s(v) &= \sup\{ |\ip{w}{v}| : w\in\mc B(E)\otimes\mc B(F),
\beta(w)\leq1\} \\ &= \sup\{ |\ip{\sigma}{w(\tau)}| : \|w\|\leq 1\}
\leq \|\sigma\| \|\tau\|\leq 1.
\end{align*}
It hence follows that the norm of $u$ as a member of the dual
space of $\mc B(E)_* \proten_{\beta^s} \mc B(F)_*$ is at least
$\beta(u)-\epsilon$.  The proof is complete, as $\epsilon>0$
was arbitrary.
\end{proof}

It would be nice if we could find a universal way to take the tensor
product of two dual Banach algebras.  For example, the projective
tensor product of two Banach algebras always gives a Banach algebra
(although it is not always the most natural norm to use, for example
for $C(K)$ spaces).  This problem is related to the fact that we
cannot find maximal or minimal WAP-crossnorms.

Let $\mc A$ and $\mc B$ be WAP algebras, and let $\alpha$ be a
WAP-crossnorm on $\mc A\otimes\mc B$.  It would be natural if there
was some connection between the DBA enveloping algebra
$\wap( (\mc A\proten_\alpha\mc B)' )'$ and the dual Banach algebra
tensor product $\wap(\mc A')' \overline\otimes_\alpha \wap(\mc B')'$.
However, this latter algebra is the dual of $\wap(\mc A') \proten_{\alpha^s}
\wap(\mc B')$, which is only a norming submodule of
$\wap( (\mc A\proten_\alpha\mc B)' )$.  Hence, in general,
$\wap(\mc A')' \overline\otimes_\alpha \wap(\mc B')'$ is only a
quotient of $\wap( (\mc A\proten_\alpha\mc B)' )'$.
Example~\ref{DBA_not_ten} shows that this is true even for
commutative C$^*$-algebras.

\subsection{Application to Connes-amenability}

\begin{theorem}\label{ten_inj}
Let $\mc A$ and $\mc B$ be Connes-amenable dual Banach algebras, and
let $\beta$ be a reasonable crossnorm on $\mc A_*\otimes\mc B_*$
which turns $(\mc A_*\proten_{\beta}\mc B_*)' =
\mc A\overline\otimes\mc B$ into a dual Banach algebra containing
$\mc A\proten_{\beta^s}\mc B$ as a weak$^*$-dense Banach algebra
(so that $\beta^s$ is an algebra crossnorm on $\mc A\otimes\mc B$).
Then $\mc A\overline\otimes\mc B$ is Connes-amenable.
\end{theorem}
\begin{proof}
Let $E$ be a reflexive Banach space, and let $\pi:\mc A\overline\otimes
\mc B\rightarrow\mc B(E)$ be a weak$^*$-continuous representation.
As in the proof of Proposition~\ref{when_wap_ca}, we claim that
$\pi(\mc A\overline\otimes\mc B)^c = \pi(\mc A\otimes\mc B)^c$, which
follows as $\pi$ is weak$^*$-continuous and $\mc A\otimes\mc B$ is
weak$^*$-dense in $\mc A\overline\otimes\mc B$.  We wish to show that
there is a quasi-expectation $\mc Q:\mc B(E) \rightarrow
\pi(\mc A\otimes\mc B)^c$.

As $\mc A$ and $\mc B$ are Connes-amenable, they are unital, with
units $e_{\mc A}$ and $e_{\mc B}$, say.  We may define a homomorphism
$\phi:\mc A\rightarrow\mc A\overline\otimes\mc B$ by
$\phi(a) = a\otimes e_{\mc B}$ for $a\in\mc A$.  Then, for
$\tau=\sum_{k=1}^n \mu_k\otimes\lambda_k\in\mc A_*\otimes\mc B_*$,
we see that
\[ \ip{\phi(a)}{\tau} = \ip{a\otimes e_{\mc B}}{\tau}
= \sum_{k=1}^n \ip{a}{\mu_k}\ip{e_{\mc B}}{\lambda_k}
= \ip{a}{\phi_*(\tau)} \qquad (a\in\mc A), \]
where $\phi_*(\tau) = \sum_{k=1}^n \ip{e_{\mc B}}{\lambda_k}\mu_k
\in\mc A_*$.  Clearly $\phi_*$ is bounded (as $\beta$ is a reasonable
crossnorm), so $\phi_*$ extends to $\mc A_*\proten_\beta\mc B_* =
(\mc A\overline\otimes\mc B)_*$, and we see that $\phi_*' = \phi$,
so that $\phi$ is weak$^*$-continuous.  A similar remark holds for $\mc B$.

Consider the representation $\pi_{\mc A}:\mc A\rightarrow\mc B(E)$ given by
$\pi_{\mc A}(a) = \pi(a\otimes e_{\mc B})$.  This is weak$^*$-continuous
by the preceding paragraph, so identify $\mc A$ with its image in $\mc B(E)$.
As $\mc A$ is Connes-amenable, there is a quasi-expectation $\mc Q_{\mc A}:
\mc B(E)\rightarrow\mc A^c$.  Analogously, there exists a quasi-expectation
$\mc Q_{\mc B}:\mc B(E)\rightarrow\mc B^c$.  Notice that
$\mc A\subseteq\mc B^c$, $\mc B\subseteq\mc A^c$ and $(\mc A\otimes\mc B)^c
= \mc A^c \cap \mc B^c$.

Let $\mc Q = \mc Q_{\mc B}\mc Q_{\mc A}$, so that $\mc Q$ is bounded,
and $\mc Q(a)=a$ for each $a\in\mc A^c \cap \mc B^c$.  Let $T\in\mc B(E)$,
let $x = \mc Q_{\mc A}(T) \in\mc A^c$, and let $b = \mc Q_{\mc B}(x)=
\mc Q(T) \in \mc B^c$.  Let $a\in\mc A$, so that $ax = xa$, and
as $\mc A\subseteq\mc B^c$,
\[ \mc Q_{\mc B}( a(b-x) ) = a \mc Q_{\mc B}(b-x) = 0
= \mc Q_{\mc B}(b-x) a = \mc Q_{\mc B}( (b-x)a ). \]
As $ab,ba\in\mc B^c$, we see that $ab = \mc Q_{\mc B}(ab) = \mc Q_{\mc B}(ba)
= ba$, so we conclude that $b\in\mc A^c$.  Thus $\mc Q$ maps into
$\mc A^c\cap\mc B^c$, and so we conclude that $\mc Q$ is a projection
onto $\mc A^c\cap\mc B^c$.
Now let $a,b\in\mc A^c\cap\mc B^c$, and let $T\in\mc B(E)$, so that
\[ \mc Q(aTb) = \mc Q_{\mc B}\mc Q_{\mc A}(aTb)
= \mc Q_{\mc B}\big( a\mc Q_{\mc A}(T) b\big)
= a \mc Q(T) b, \]
so we conclude that $\mc Q$ is a quasi-expectation, as required.
\end{proof}

Notice that this proof will also show that $\mc C$ is Connes-amenable,
whenever $\mc C$ is a dual Banach algebra containing
$\mc A\otimes\mc B$ as a dense subalgebra, and is such that
the map $\mc A\rightarrow\mc C; a\mapsto a\otimes
e_{\mc B}$ is weak$^*$-continuous (and similarly for $\mc B$).

\begin{corollary}
Let $r,s\in (1,\infty)$, and let $\alpha$ be some quasi-uniform
crossnorm on $\ell^r\otimes\ell^s$ such that $\ell^r\proten_\alpha\ell^s$
is reflexive, and $(\ell^r)'\otimes(\ell^s)'$ is dense in
$(\ell^r\proten_\alpha\ell^s)'$.  Then $\mc B(\ell^r
\proten_\alpha\ell^s)$ is Connes-amenable.
\end{corollary}
\begin{proof}
By Proposition~\ref{for_op_spaces}, we have that
$\mc B(\ell^r\proten_\alpha\ell^s) = \mc B(\ell^r)
\overline\otimes_\alpha\mc B(\ell^s)$.  Then the theorem applies,
as $\mc B(\ell^r)$ and $\mc B(\ell^s)$ are Connes-amenable, by
results in \cite{Runde1}.
\end{proof}

This corollary is comparable to \cite[Theorem~2.2]{GJW},
as a quasi-uniform tensor norm is \emph{tight} in the sense of
\cite{GJW}, and by \cite{Runde1}, the amenability of $\mc A(E)$
is equivalent to the Connes-amenability of $\mc B(E)$, at least
when $E$ is reflexive and has the approximation property.

It is interesting to note that our proof of Theorem~\ref{ten_inj}
is rather more algebraic than an analogous von Neumann result
(compare to \cite[Chapter~XV, Proposition~3.2]{tak3}).
Our approach is more in line with that of \cite[Proposition~6.3.17]{RundeBook}.

\section{Acknowledgements}

The author would like to thank the anonymous reviewer for much
helpful advice.

\end{document}